\documentclass[10pt,a4paper]{amsart}
\usepackage{amsmath, graphicx, subfigure, epsfig,amssymb}
\usepackage{xcolor}
\usepackage[normalem]{ulem}
\numberwithin{equation}{section}

\newcommand{\e}{\epsilon}
\newcommand{\ga}{\gamma}
\newcommand{\de}{\delta}
\newcommand{\br}{\mathbb{R}}

\newcommand{\ik}{\varphi}
\newcommand{\pa}{\partial}
\newcommand{\bt}{\beta}
\newcommand{\al}{\alpha}
\newcommand{\la}{\lambda}

\newcommand{\be}{\begin{equation}}
\newcommand{\ee}{\end{equation}}

\newcommand{\tth}{\tilde\theta}
\newcommand{\tal}{\tilde\al}
\newcommand{\hal}{\hat\al}

\newcommand{\dd}{\text{d}}
\newcommand{\frec}{f_\e^{\text{rec}}}

\newcommand{\CA}{\mathcal{A}}

\newcommand{\s}{\mathcal S}

\newcommand{\R}{\mathcal R}

\newcommand{\CH}{\mathcal H}
\newcommand{\BZ}{\mathbb Z}
\newcommand{\Iz}{I^{(0)}}

\newcommand{\DTB}{\text{DTB}}

\newtheorem{theorem}{Theorem}
\newtheorem{remark}{Remark}
\newtheorem{lemma}{Lemma}
\newtheorem{corollary}{Corollary}
\newtheorem{definition}{Definition}

\begin{document}

\title[Analysis of resolution]{\textcolor[rgb]{0.9, 0, 0}{Novel resolution analysis for\\ the Radon transform in $\mathbb R^2$ for\\ functions with rough edges}}
\author[A Katsevich]{Alexander Katsevich$^1$}
\thanks{$^1$This work was supported in part by NSF grant DMS-1906361. Department of Mathematics, University of Central Florida, Orlando, FL 32816 (Alexander.Katsevich@ucf.edu). }

\begin{abstract} 
Let $f$ be a function in $\mathbb R^2$, which has a jump across a smooth curve $\mathcal S$ with nonzero curvature. We consider a family of functions $f_\epsilon$ with jumps across a family of curves $\mathcal S_\epsilon$. Each $\mathcal S_\epsilon$ is an $O(\epsilon)$-size perturbation of $\mathcal S$, which scales like $O(\epsilon^{-1/2})$ along $\mathcal S$. Let $f_\epsilon^{\text{rec}}$ be the reconstruction of $f_\epsilon$ from its discrete Radon transform data, where $\epsilon$ is the data sampling rate. 
A simple asymptotic (as $\epsilon\to0$) formula to approximate $f_\epsilon^{\text{rec}}$ in any $O(\epsilon)$-size neighborhood of $\mathcal S$ was derived heuristically in an earlier paper of the author. Numerical experiments revealed that the formula is highly accurate even for nonsmooth (i.e., only H{\"o}lder continuous) $\mathcal S_\epsilon$. In this paper we provide a full proof of this result, which says that the magnitude of the error between $f_\epsilon^{\text{rec}}$ and its approximation is $O(\epsilon^{1/2}\ln(1/\epsilon))$. The main assumption is that the level sets of the function $H_0(\cdot,\epsilon)$, which parametrizes the perturbation $\mathcal S\to\mathcal S_\epsilon$, are not too dense. 
\end{abstract}
\maketitle

\section{Introduction}\label{sec_intro}

\subsection{Local resolution analysis: original and new} Let $f$ be an unknown function in $\br^2$, and $\s$ be some surface. We assume that $f$ has a jump discontinuity across $\s$, and $f$ is sufficiently smooth otherwise. Let $\frec$ be a reconstruction from discrete tomographic data (i.e., the classical Radon transform of $f$), where $\e$ represents the data sampling rate. The reconstruction is computed by substituting interpolated data into a ``continuous'' inversion formula. In many applications it is important to know the resolution of reconstruction from discrete data, including medical imaging, materials science, and nondestructive testing.

In \cite{Katsevich2017a} the author initiated analysis of resolution,  called {\it local resolution analysis}, by focusing specifically on the behavior of $\frec$ near $\s$. One of the main results of \cite{Katsevich2017a} is the computation of the limit 
\be\label{tr-beh}
\DTB(\check x)=\lim_{\e\to0}\frec(x_0+\e\check x)
\ee
in a 2D setting under the assumptions that (a) $\s$ a sufficiently smooth curve with nonzero curvature, (b) $f$ has a jump discontinuity across $\s$, (c) $x_0\in\s$ is generic, and (d) $\check x$ is confined to a bounded set. It is important to emphasize that both the size of the neighborhood around $x_0$ and the data sampling rate go to zero simultaneously in \eqref{tr-beh}. The limiting function $\DTB(\check x)$, which we call the discrete transition behavior (or DTB for short), contains complete information about the resolution of reconstruction. The practical use of the DTB is based on the relation
\be\label{DTB orig use}
\frec(x_0+\e\check x)=\DTB(\check x)+\text{error term}.
\ee
When $\e>0$ is sufficiently small, the error term is negligible, and $\DTB(\check x)$, which is given by a simple formula, is an accurate approximation to $\frec$. Numerical experiments reported in \cite{Katsevich2017a} demonstrate that the error term in \eqref{DTB orig use} is indeed quite small for realistic values of $\e$. Local resolution analysis was extended to much more general settings in subsequent papers \cite{kat19a, kat20b, kat20a, Kats_21resol}. 

The functions, which have been investigated in the framework of local resolution analysis are, for the most part, nonsmooth across smooth surfaces. On the other hand, in many applications discontinuities of $f$ occur across non-smooth (rough) surfaces. Examples include soil and rock imaging, where the surfaces of cracks and pores and boundaries between adjacent regions with different properties are highly irregular and frequently simulated by fractals \cite{Anovitz2015, GouyetRosso1996, Li2019, soilfractals2000, PowerTullis1991, Renard2004, Zhu2019}. 

It was proven in \cite{Katsevich2021b} that the original approach to resolution analysis based on \eqref{tr-beh} still works for functions with jumps across nonsmooth curves
(i.e., H{\"o}lder continuous with some exponent $\ga\in(0,1]$). 
More precisely, we consider a family of functions $f_\e$, which have a jump across $\s_\e$, where $\s_\e$ is a perturbation of a smooth and convex curve $\s$, and $\s_\e$ is not necessarily smooth, and then show that under certain assumptions on $f_\e$ and $\s_\e$, the error term in \eqref{DTB orig use} goes to zero as $\e\to0$. 
Nevertheless, numerical experiments reported in \cite{Katsevich2021b} demonstrated that this approach is not entirely satisfactory: a significant mismatch between the reconstruction $\frec(x_0+\e\check x)$ and its approximation $\DTB(\check x)$ was observed. This means that the error term in \eqref{DTB orig use} decays slowly as $\e\to0$ when $\s_\e$ is rough. 

To overcome this problem, a new approach to local resolution analysis was developed in \cite{Katsevich2021b}. This approach is based on allowing the DTB to depend on $\e$:
\be\label{DTB new use}
\frec(x_0+\e\check x)=\DTB_{new}(\check x,\e)+\text{error term}.
\ee
The idea is that with the new DTB being more flexible (due to its $\e$-dependence), the error term in \eqref{DTB new use} can be smaller than the one in \eqref{DTB orig use}. The new DTB proposed in \cite{Katsevich2021b} (see eq. \eqref{via_kernel_lim} below) is given by the convolution of an explicitly computed and suitably scaled kernel with $f$. Thus, analysis of resolution based on \eqref{DTB new use} is as simple as the one based on \eqref{DTB orig use}, and it can be used quite easily to investigate partial volume effects and resolution in the case of rough boundaries. Numerical experiments with the new DTB presented in \cite{Katsevich2021b} show an excellent match between $\DTB_{new}$ and actual reconstructions even when $\s_\e$ is fractal.  

To prove that the new DTB works well for nonsmooth $\s_\e$, it is not sufficient to show that the error term in \eqref{DTB new use} goes to zero as $\e\to0$. We need to establish that the magnitude of the error is {\it independent} of how rough $\s_\e$ is, i.e. independent of its H{\"o}lder exponent $\ga$. This proved to be a difficult task. In \cite{Katsevich2021b} it was conjectured that the error term in \eqref{DTB new use} is $O(\e^{1/2}\ln(1/\e))$, and a partial results towards proving this conjecture was established. In this paper we provide a full proof of the conjecture. 

Our approach is asymptotic in nature. As mentioned above, instead of considering a single, fixed function $f$, we consider a family of functions $f_\e$. Each $f_\e$ has a jump discontinuity across $\s_\e$, where $\s_\e$ is a small and suitably scaled perturbation of a smooth and convex curve $\s$. The parametrization of $\s_\e$ is given by $\theta\to y(\theta)+\e H_0(\e^{-1/2}\theta,\e)\vec\theta$, where $\vec\theta=(\cos\theta,\sin\theta)$, and $\theta$ runs over an interval. Here $\theta\to y(\theta)\in\s$ is the parameterization of $\s$ such that $\vec\theta$ is orthogonal to $\s$ at $y(\theta)$. The main assumption is that the function $H_0$ (more precisely, the family of functions $H_0(\cdot,\e)$), which determines the perturbation, has level sets that are not too dense. This assumption does allow fairly nonsmooth $H_0$. For example, in \cite{Katsevich2021b} we construct a function on $\br$, whose level sets are not too dense as required, which is H{\"o}lder continuous with exponent $\ga$ for any prescribed $0<\ga<1$, but which is not H{\"o}lder continuous with any exponent $\ga'>\ga$ on a dense subset of $\br$. 
Our construction ensures that the size of the perturbation is $O(\e)$ in the direction normal to $\s$, and the perturbation scales like $\e^{-1/2}$ in the direction tangent to $\s$. 

The ideas behind the proof in this paper are quite different from those used in the original approach \cite{Katsevich2017a, kat19a, kat20b, kat20a, Kats_21resol}. The latter proofs revolve around the smoothness of the singular support of $f$. The new proofs are based on the phenomenon of destructive interference, which is described mathematically as the smallness of certain exponential sums. The assumption about the level sets of $H_0$ is what allows for the destructive interference to occur.

\subsection{Practical application of our results} An important application of our results is in micro-CT (i.e., CT capable of achieving micrometer resolution), which is a valuable tool for imaging of rock samples extracted from wells. The reconstructed images can then be used to investigate properties of the samples.
A collection of numerical methods that determine various rock properties using digital cores is collectively called Digital Rock Physics (or, DRP) \cite{Fredrich2014, Saxena2019}. Here the term ``digital core'' refers to a digital representation of the rock sample (rock core) obtained as a result of micro-CT scanning, reconstruction, and image analysis (segmentation and classification, feature extraction, etc.) \cite{Guntoro2019}. DRP ``is a rapidly advancing technology that relies on digital images of rocks to simulate multiphysics at the pore-scale and predict properties of complex rocks (e.g., porosity, permeability, compressibility). ... For the energy industry, DRP aims to achieve more, cheaper, and faster results as compared to conventional laboratory measurements.'' \cite{Saxena2019}. Furthermore, ``The simulation of various rock properties based on three-dimensional digital cores plays an increasingly important role in oil and gas exploration and development. {\it The accuracy of 3D digital core reconstruction is important for determining rock properties}.''   \cite{Zhu2019} (italic font is added here).

As stated above, boundaries between regions with different properties inside the rock (e.g., between the solid matrix and the pore space) are typically rough (see also \cite{Cherk2000}), i.e., they contain features across a wide range of scales, including scales below what is accessible with micro-CT. 
Given that the quality of micro-CT images is critical for accurate DRP, effects that degrade the resolution of micro-CT (e.g., due to finite data sampling) and how these effects manifest themselves in the presence of rough boundaries require careful investigation. Once fully understood and quantified, these effects can be accounted for at the step of image analysis, thereby leading to more accurate DRP results.

\subsection{Related results. Organization of the paper} To situate our paper in a more general context, not much is known about how the Radon transform acts on distributions with complicated singularities. A recent literature search reveals a small number of works, which investigate the Radon transform acting on random fields \cite{JainAnsary1984, Sanz1988, Medina2020}. The author did not find any publication on the Radon transform of characteristic functions of domains with rough boundaries. This appears to be the first paper that contains a result on the Radon transform of functions with rough edges. 

An alternative way to study resolution of tomographic reconstruction is based on sampling theory. Applications of the classical sampling theory to Radon inversion are in papers such as \cite{nat93, pal95, far04}, just to name a few. Analysis of sampling for distributions with semiclasical singularities is in \cite{stef20, Monard2021}. This line of work determines the sampling rate required to reliably recover features of a given size and describes aliasing artifacts if the sampling requirements are violated.

The paper is organized as follows. In section~\ref{sec:a-prelims} we describe the setting of the problem, state the definition of a generic point, formulate all the assumptions (including assumptions about the perturbation $H_0$), and formulate the main result (Theorem~\ref{main-res}). The beginning of the proof is in section~\ref{sec:beg proof}. We consider three cases: 
\begin{itemize}
\item[(A)] $x_0\in\s$; 
\item[(B)] $x_0\not\in\s$, and there is a line through $x_0$, which is tangent to $\s$; and 
\item[(C)] $x_0\not\in\s$, and no line through $x_0$ is tangent to $\s$. 
\end{itemize}
Sections~\ref{sec:prep}--\ref{sec: sum_k I} contain a nearly complete proof of the theorem in case (A). What is left is one additional assertion, which is proven in a later section. Likewise, Sections~\ref{sec: rem prep} and \ref{sec: sum_k II} contain a nearly complete proof of the theorem in case (B). Originally, case (C) is proven in \cite{Katsevich2021b}. Its proof is given in appendix~\ref{sec new ker prf} to make the paper self-contained. The final assertion of the theorem in cases (A) and (B) is proven in section~\ref{m=0 case}. Proofs of all lemmas and some auxiliary results are in appendices~\ref{sec:loc sing proofs}--\ref{sec: spec cases}.

%

\section{Preliminaries}\label{sec:a-prelims}

\begin{figure}[h]
{\centerline
{\epsfig{file=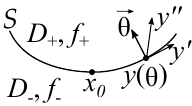, width=4.0cm}}}
\caption{Illustration of a function $f$ discontinuous across $\s$.}
\label{fig:jump}
\end{figure}

Consider a compactly supported function $f(x)$ in the plane, $x\in\br^2$. Set $\s:=\{x\in\br^2:f\not\in C^2(U)\text{ for any open }U\ni x\}$. We suppose that 
\begin{itemize}
\item[f1.] For each $x_0\in\s$ there exist a neighborhood $U\ni x_0$, domains $D_\pm$, and functions $f_\pm\in C^2(\br^2)$ such that
\begin{equation}\label{f_def}\begin{split}
& f(x)=\chi_{D_-}(x) f_-(x)+\chi_{D_+}(x) f_+(x),\ x\in U\setminus \s,\\
& D_-\cap D_+=\varnothing,\ D_-\cup D_+=U\setminus \s,
\end{split}
\end{equation}
where $\chi_{D_\pm}$ are the characteristic functions of $D_\pm$,  
\item[f2.] $\s$ is a $C^4$ curve.
\end{itemize}
Conditions f1, f2 describe a typical function, which has a jump discontinuity across a smooth curve (see Figure~\ref{fig:jump}).

The discrete tomographic data are given by
\be\label{data_eps}
\hat f_\e(\al_k,p_j):=\frac1\e \iint w\left(\frac{p_j-\vec\al_k\cdot y}{\e}\right)f(y)\dd y,\
p_j=\bar p+j\Delta p,\ \al_k=\bar\al+k\Delta\al,
\ee
where $w$ is the detector aperture function, $\Delta p=\e$, $\Delta\al=\kappa\e$, and $\kappa>0$, $\bar p$, $\bar\al$, are fixed. 
Here and below, $\vec \al$ and $\al$ in the same equation are always related by $\vec\al=(\cos\al,\sin\al)$. The same applies to $\vec\theta=(\cos\theta,\sin\theta)$ and $\theta$. Sometimes we also use $\vec\theta^\perp=(-\sin\theta,\cos\theta)$.

\noindent{\bf Assumptions about the Aperture Function $w$:}
\begin{itemize}
\item[AF1.] $w$ is even and $w\in C_0^{\lceil \bt \rceil+1}(\br)$ (i.e., $w$ is compactly supported, and $w^{(\lceil \bt\rceil+1)}\in L^\infty(\br)$) for some $\bt\ge 3$; and 
\item[AF2.] $\int w(p)dp=1$.
\end{itemize}
Here $\lceil \bt\rceil$ is the ceiling function, i.e. the integer $n$ such that $n-1<\bt\le n$. The required value of $\bt$ is stated below in Theorem~\ref{main-res}. Later we also use the floor function $\lfloor \bt\rfloor$, which gives the integer $n$ such that $n\le\bt<n+1$, and the fractional part function $\{\bt\}:=\bt-\lfloor \bt\rfloor$. 

Reconstruction from discrete data is achieved by the formula
\be\label{recon-orig}
\frec(x)=-\frac{\Delta\al}{2\pi}\sum_{|\al_k|\le \pi/2} \frac1\pi \int \frac{\pa_p\sum_j\ik\left(\frac{p-p_j}\e\right)\hat f_\e(\al_k,p_j)}{p-\al_k\cdot x}\dd p,
\ee
where $\ik$ is an interpolation kernel.

\noindent{\bf Assumptions about the Interpolation Kernel $\ik$:}
\begin{itemize}
\item[IK1.] $\ik$ is even, compactly supported, and its Fourier transform satisfies $\tilde\ik(\la)=O(|\la|^{-(\bt+1)})$, $\la\to\infty$;
\item[IK2.] $\ik$ is exact up to order $1$, i.e. 
\be\label{exactness}
\sum_{j\in\mathbb Z} j^m\ik(u-j)\equiv u^m,\ m=0,1,\ u\in\br.
\ee
\end{itemize}
Here $\bt$ is the same as in assumption AF1. As is easily seen, assumption IK2 implies $\int\ik(p)dp=1$.
See Section IV.D in \cite{btu2003}, which shows that $\ik$ with the desired properties can be found for any $\bt>0$ (i.e., for any regularity of $\ik$).


We can assume that the coordinates are selected so that $\vec\theta_0:=(1,0)$. Suppose $\s$ is parametrized by $[-a,a]\ni\theta\to y(\theta)\in\s$. Here $y(\theta)$ is the point such that $\vec\theta\cdot y'(\theta)=0$. Suppose $\vec\theta^\perp\cdot y'(\theta)<0$ and $\vec\theta$ points towards the center of curvature of $\s$. Thus, $\vec\theta^\perp=-y'(\theta)/|y'(\theta)|$ and $\R(\theta)=\vec\theta\cdot y''(\theta)>0$, $|\theta|\le a$, where $\R(\theta)$ is the radius of curvature of $\s$ at $y(\theta)$ (see Figure~\ref{fig:jump}). If $x_0\in\s$, we assume $x_0=y(0)$.

%

Let $H_0(s,\e)$, $s\in\br$, be a family of functions defined for all $\e>0$ sufficiently small. 
Define the function
\be\label{chi fn}
\chi(t,r):=\begin{cases} 1,&0< t\le r,\\
-1,&r\le t < 0,\\
0,&\text{otherwise}.\end{cases}
\ee
Suppose $H_0$ has the following properties.  

There exist constants $c,\rho,L>0$, which are independent of $\e$, such that for all $\e>0$ sufficiently small:
\begin{itemize}
\item[$H1$.] $|H_0(u,\e)|\le c$ for all $u\in\br$;
\item[$H2$.] The function $(t,u)\to \chi(t,H_0(u,\e))$ is measurable in $\br^2$; and
\item[$H3$.] For any interval $I$ of any length $L\ge L_0$ the set 
\be\label{U intervals}
\{u\in I: \text{sgn}(t)(H_0(u,\e)-t)\ge 0\}
\ee
is either empty or a union of no more than $\rho L$ intervals $U_n$, $\text{dist}(U_{n_1},U_{n_2})>0$, $n_1\not=n_2$, for almost all $t\not=0$. 
\end{itemize}
Our assumptions allow $H_0$ to be discontinuous. The endpoints of $U_n$ are denoted $u_{2n}$ and $u_{2n+1}$: $\overline{U}_n=[u_{2n},u_{2n+1}]$, where the bar denotes closure. The distance between the intervals is positve, so $u_n<u_{n+1}$ for all $n$. The intervals $U_n$ and the points $u_n$ depend on $t$, $I$, and $\e$.
If $H_0$ is continuous, then for each $t$ and $\e>0$ the collection of $u_n$'s is simply the level set $\{u\in\br:H_0(u,\e)=t\}$. Therefore, condition $H3$ can be informally interpreted as saying that the level sets of $H_0$ are not too dense, i.e. they do not become infinitely dense on intervals. 

Assumptions $H1$ and $H2$ imply that $\int_I\int_{\br}|\chi(t,H_0(u,\e))|\dd t\dd u$ is well-defined and bounded for any bounded interval $I$. Hence, by the Fubini theorem and $H3$,
\be\label{switch}\begin{split}
\int_I\int_{\br} g(t,u)\chi(t,H_0(u,\e))\dd t\dd u&=\int_{\br}\int_I g(t,u)\chi(t,H_0(u,\e))\dd u\dd t\\
&=\int_{\br}\text{sgn}(t)\sum_n\int_{U_n} g(t,u)\dd u\dd t
\end{split}
\ee
for any sufficiently regular function $g$.

\begin{figure}[h]
{\centerline
{\epsfig{file=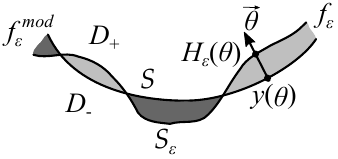, width=4.5cm}}}
\caption{Illustration of the perturbation $\s\to\s_\e$ and the function $f_\e$, which is supported in the shaded regions.}
\label{fig:perturbation}
\end{figure}

In what follows, the dependence of $H_0$ on $\e$ is omitted from notation for simplicity. Set $H_\e(\theta):=\e H_0(\e^{-1/2}\theta)$ and define also
\be\label{main-fn}
f_\e(x)=(f_+(x)-f_-(x))\chi(t,H_\e(\theta)),\ x=y(\theta)+t\vec\theta.
\ee
As is easily seen, $f_\e^{mod}(x):=f(x)-f_\e(x)$ is the function, in which $\s$ is modified by $H_\e$, and $f_\e$ is the corresponding perturbation, see Figure~\ref{fig:perturbation}. At the points where $H_\e(\theta)>0$, a small region is removed from $D_+$ and added to $D_-$ (see lighter shaded regions in Figure~\ref{fig:perturbation}). At the points where $H_\e(\theta)<0$, a small region is removed from $D_-$ and added to $D_+$ (see darker shaded regions in Figure~\ref{fig:perturbation}). The magnitude of the perturbation is $O(\e)$. Let $\s_\e$ denote the perturbed boundary. Thus, $f_\e^{mod}(x)$ is discontinuous across $\s_\e$ instead of $\s$.

In \cite{Katsevich2017a} we obtained the DTB in the case of a sufficiently smooth $\s$. By linearity, we can ignore the original function $f$ and consider the reconstruction of only the perturbation $f_\e$. 

Let $\frec$ denote the reconstrution of only the perturbation $f_\e$ \eqref{main-fn}. By \eqref{data_eps}, \eqref{recon-orig},
\be\label{via_kernel}
\frec(x)=-\frac{\Delta\al}{2\pi}\frac1{\e^2}\sum_{|\al_k|\le \pi/2}\sum_j \CH\ik'\left(\frac{\vec\al_k\cdot x-p_j}\e\right) \iint w\left(\frac{p_j-\vec\al_k\cdot y}{\e}\right)f_\e(y)\dd y.
\ee
Following \cite{Katsevich2021b}, replace the sums with respect to $k$ and $j$ with integrals:
\be\label{via_kernel_lim}\begin{split}
&\frec(x)\approx\frac{1}{\e^2}\iint K\left(\frac{x-y}\e\right)f_\e(y)\dd y,\
K(z):=-\frac{1}{2\pi}\int_0^\pi (\CH\ik'*w)(\vec\al\cdot z) \dd\al.
\end{split}
\ee
As is easily seen, $K$ is radial and compactly supported. 

For a real number $s$, let $\langle s\rangle$ denote the distance from $s$ to the nearest integer:
$\langle s\rangle:=\min_{l\in\mathbb Z}|s-l|=\min(\{s\},1-\{s\})$. The following definition is in \cite[p. 121]{KN_06} (after a slight modification in the spirit of \cite[p. 172]{Naito2004}).

\begin{definition} Let $\eta>0$. The irrational number $s$ is said to be of type $\eta$ if for any $\eta_1>\eta$, there exists $c(s,\eta_1)>0$ such that
\be\label{type ineq}
m^{\eta_1}\langle ms\rangle \ge c(s,\eta_1) \text{ for any } m\in\mathbb N.
\ee
The irrational number $s$ is said to be of constant type or badly approximable if there exists $c(s)>0$ such that
$m\langle ms\rangle \ge c(s)$ for any $m\in\mathbb N$.
\end{definition}
See also \cite{Naito2004}, where the numbers which satisfy \eqref{type ineq} are called $(\eta-1)$-order Roth numbers. It is known that $\eta\ge1$ for any irrational $s$. The set of irrationals of each type $\eta \ge 1$ is of full measure in the Lebesgue sense, while the class of constant type is dense in the space of real numbers, but it is of null measure \cite{Naito2004}.

In the rest of the paper we consider points $x_0$ such that 
\begin{itemize}
\item[P1.] No line, which is tangent to $\s$ at a point where the curvature of $\s$ is zero, passes through $x_0$;
\item[P2.] The line through the origin and $x_0$ is not tangent to $\s$; 
\item[P3.] $\kappa|x_0|$ is irrational and of finite type; and
\item[P4.] If $x_0\in\s$, i.e. $x_0=y(0)$, then the number $\kappa \theta_0^\perp\cdot x_0$ is irrational and of finite type.
\end{itemize}

Assumption P2 is well-defined, since the geometry of the tomographic data determines the origin. By P4, $x_0\not=0$. Clearly, the set of such $x_0$ is dense in the plane. Let $\eta_0$ denote the type of $\kappa|x_0|$ if $x_0\not\in\s$, and the larger of the two types (in P3 and P4) -- if $x_0\in\s$. Our main result is the following theorem. 

\begin{theorem}\label{main-res} Suppose 
\begin{enumerate}
\item $f$ satisfies conditions f1, f2; 
\item the detector aperture function $w$ satisfies conditions AF1, AF2; 
\item the interpolation kernel $\ik$ satisfies conditions IK1, IK2; 
\item the perturbation $H_0$ satisfies conditions $H1$--$H3$; and 
\item point $x_0$ satisfies conditions P1--P4. 
\end{enumerate}
If $\bt>\eta_0+2$, one has:
\be\label{full conj}
\frec(x_0+\e\check x)=\frac1{\e^2}\iint K\left(\frac{(x_0+\e\check x)-y}\e\right)f_\e(y)\dd y+O(\e^{1/2}\ln(1/\e)),\ \e\to0,
\ee 
where $K$ is given by \eqref{via_kernel_lim}, and the big-$O$ term is uniform with respect to $\check x$ in any compact set.
\end{theorem}


\section{Beginning of the proof of Theorem~\ref{main-res}}\label{sec:beg proof}

By linearity, in what follows we can consider only one domain $U$ and make the following assumptions: 
\begin{enumerate}
\item $\text{supp}(f)\subset U$,  
\item $\s$ is sufficiently short;
\item $f\equiv0$ in a neighborhood of the endpoints of $\s$.
\end{enumerate}
By assumption (3) above, $f(x)\equiv0$ in a neighborhood of $y(\pm a)$. The requirements on the smallness of $a$ are formulated later as needed.

Throughout the paper we use the following convention. If an inequality involves an unspecified constant $c$, this means that the inequality holds for some $c>0$. If an inequality (or, a string of inequalities) involves multiple unspecified constants $c$, then the values of $c$ in different places can be different. If some additional information about the value of $c$ is necessary (e.g., $c\gg1$ or $c>0$ small), then it is stated.

Following \cite{Katsevich2021b}, consider the function (see \eqref{via_kernel})
\be\label{recon-ker-0}
\psi(q,t):=\sum_j (\CH\ik')(q-j)w(j-q-t).
\ee
Then
\be\label{psi-props-0}\begin{split}
&\psi(q,t)=\psi(q+1,t),\  q,t\in\br;\quad
\psi(q,t)=O(t^{-2}),\ t\to\infty,\ q\in\br;\\
&\int \psi(q,t)\dd t\equiv 0,\ q\in\br.
\end{split}
\ee
The last property follows from IK2 (see \eqref{exactness}).
By \eqref{psi-props-0}, we can represent $\psi$ in terms of its Fourier series:
\be\label{four-ser}\begin{split}
\psi(q,t)&=\sum_m \tilde\psi_m(t) e(-mq),\ e(q):=\exp(2\pi i q),\\
\tilde\psi_m(t) &=\int_0^1 \psi(q,t)e(mq)\dd q=\int_\br (\CH\ik')(q) w(-q-t)e(mq)\dd q.
\end{split}
\ee
Introduce the function $\rho(s):=(1+|s|^\bt)^{-1}$, $s\in\br$. 
\begin{lemma}\label{lem:psi psider}
One has
\be\label{four-coef-bnd}
|\tilde\psi_m(t)|,|\tilde\psi_m'(t)| \le c\rho(m)(1+t^2)^{-1}.
\ee
\end{lemma}
By the lemma, the Fourier series for $\psi$ converges absolutely.

From \eqref{data_eps}, \eqref{via_kernel}, \eqref{recon-ker-0},  and \eqref{four-ser}, the reconstructed image becomes
\be\label{recon-ker-v2}
\begin{split}
&\frec(x)
=-\frac{\Delta\al}{2\pi}\sum_m \sum_{|\al_k|\le \pi/2} e\left(-m \frac{\vec\al_k\cdot x-\bar p}\e\right)A_m(\al_k,\e), \\
&A_m(\al,\e):=\e^{-2}\iint \tilde\psi_m\left(\frac{\vec\al\cdot (y-x)}\e\right)f_\e(y)\dd y.
\end{split}
\ee
To obtain \eqref{via_kernel_lim}, in \eqref{via_kernel} we should be able to replace the sum with respect to $k$ by an integral with respect to $\al$ and ignore all $m\not=0$ terms (that make up $\psi$). We will show that 
\begin{align}\label{extra terms I}
&\Delta\al\sum_{m\not=0}\left| \sum_{|\al_k|\le\pi/2} e\left(-m \frac{\vec\al_k\cdot x}\e\right)A_m(\al_k,\e)\right|=O(\e^{1/2}\ln(1/\e)),\  \\
&\label{extra terms II}
\sum_{|\al_k|\le \pi/2}
\int_{\al_k-\Delta\al/2}^{\al_k+\Delta\al/2}  \left|A_0(\al,\e)-A_0(\al_k,\e)\right| \dd\al =O(\e^{1/2}\ln(1/\e)),
\end{align}
where $x=x_0+\e\check x$. The factor $e(m\bar p/\e)$ is dropped, because it is independent of $k$.

%

All the estimates below are uniform with respect to $\check x$, so the $\check x$-dependence of various quantities is frequently omitted from notation. Transform the expression for $A_m$ (cf. \eqref{recon-ker-v2}) by changing variables $y\to (\theta,t)$, where $y=y(\theta)+t\vec\theta$:
\be\label{A-simpl}
\begin{split}
A_m(\al,\e)=&\frac1{\e^2}\int_{-a}^a\int_0^{H_\e(\theta)}\tilde\psi_m\left(\frac{\vec\al\cdot (y(\theta)-x_0)}\e+h(\theta,\al)\right)F(\theta,t)\dd t\dd\theta,\\
F(\theta,t):=&\Delta f(y(\theta)+t\vec\theta) (\R(\theta)-t),\ h(\theta,\al):=-\al\cdot\check x+\hat t\cos(\theta-\al),
\end{split}
\ee
where $\R(\theta)-t=\text{det}(\dd y/\dd (\theta,t))>0$. Recall that $\R(\theta)$ is the radius of curvature of $\s$ at $y(\theta)$. The dependence of $h$ on $\hat t$ and $\check x$ is irrelevant and omitted from notation. 

Consider the function
\be\label{bigR1 orig}
R_1(\theta,\al):=\vec\al\cdot (y(\theta)-x_0).
\ee
Change variables $\theta=\e^{1/2}\tth$ and $t=\e\hat t$ in \eqref{A-simpl} and then use \eqref{switch}:
\be\label{A simpl orig}
\begin{split}
\e^{1/2}A_m(\al,\e)=&\int_{-a\e^{-1/2}}^{a\e^{-1/2}}\int_{\br}\tilde\psi_m\left(\e^{-1}R_1(\theta,\al)+h(\theta,\al)\right)
F(\theta,\e \hat t)\chi(\hat t,H_0(\tilde\theta))\dd\hat t\dd\tilde\theta\\
=&\int_{\br}\text{sgn}(\hat t)\sum_n\int_{U_n}\tilde\psi_m\left(\e^{-1}R_1(\theta,\al)+h(\theta,\al)\right)F(\theta,\e \hat t)\dd\tilde\theta\dd\hat t,
\end{split}
\ee
where $\theta$ is a function of $\tth$: $\theta=\e^{1/2}\tth$, and $\chi$ is defined in \eqref{chi fn}. Recall that $U_n=[u_{2n},u_{2n+1}]$ (see $H3$, where these intervals are introduced). Here and in what follows, the interval $I$ used in the construction of $U_n$'s is always $I=[-a\e^{-1/2},a\e^{-1/2}]$. The intervals $U_n$ can be closed, open, and half-closed. Since what kind they are is irrelevant, with some abuse of notation we write them is if they are closed.
The dependence of $U_n$ and $u_n$ on $\hat t$ and $\e$ is omitted from notation for simplicity. 


In view of \eqref{recon-ker-v2}, \eqref{extra terms I}, and \eqref{A simpl orig} we need to estimate the quantity
\be\label{W def orig}
\begin{split}
&W_{m}(\hat t):=\e^{1/2} \sum_{\Delta\al|k|\le\pi/2 }e\left(-m q_k\right) \left[g(\tal_k;\cdot) e\left(-m \vec\al_k\cdot \check x\right)\right],\ q_k:=\frac{\vec\al_k\cdot x_0}\e,\\
&g(\tal;m,\e,\hat t,\check x):=\sum_n\int_{U_n}\tilde\psi_m\left(\e^{-1}R_1(\theta,\al)+h(\theta,\al)\right)F(\theta,\e \hat t)\dd\tilde\theta.
\end{split}
\ee
For convenience, we express $g$ as a function of $\tal$ rather than $\al$. The sum in \eqref{extra terms I} is bounded by $\int \sum_{m\not=0}|W_{m}(\hat t)|\dd\hat t$. By property $H1$, $H_0$ is bounded, so the integral with respect to $\hat t$ is over a bounded set.

Throughout the paper we frequently used rescaled variables $\tth,\hat\theta$ and $\tal,\hal$:
\be\label{scales}
\tal=\e^{-1/2}\al,\ \hal=\al/\Delta\al,\ \tth=\e^{-1/2}\theta,\ \hat\theta=\theta/\Delta\al.
\ee
Whenever an original variable (e.g., $\al$) is used together with its rescaled counterpart (e.g., $\hal$ or $\tal$) in the same equation or sentence, they are always assumed to be related according to \eqref{scales}. The same applies to $\theta$ and its rescaled versions. The only exception is appendix~\ref{sec new ker prf}, where the relationship between $\theta$ and $\tth$ is slightly different from the one in \eqref{scales}.

We distinguish two cases: $x_0\in\s$ and $x_0\not\in\s$. In the former case, $x_0=y(0)$.
The proof of \eqref{extra terms I} is much more difficult than the proof of \eqref{extra terms II}, so we discuss the intuition behind the former. Due to the remark in the paragraph following \eqref{W def orig}, we just have to show that $\sum_{m\not=0}|W_m(\hat t)|$ satisfies the same estimate as in \eqref{extra terms I}.
Summation with respect to $m$ does not bring any complications, so we consider the sum with respect to $k$ for a fixed $m\not=0$. The factor $g(\tal_k;\cdot)$ is bounded and goes to zero as $\tal\to\infty$. In addition, the factor $e(-m q_k)$ oscillates rapidly. If the product $g(\tal_k;\cdot) e\left(-m \vec\al_k\cdot \check x\right)$ changes slowly (in the appropriate sense),  destructive interference makes the sum small. 

An additional phenomenon is that destructive interference does not work near the points $\al_{j,m}$, where the derivative of the phase $m\vec\al^\perp(\al_{j,m})\cdot x_0=j\in\BZ$. The points $\al_{j,m}$ need to be investigated separately. If $\al_{j,m}=0$ for some $j$ and $m$, \eqref{extra terms I} may fail because in this case $g(\tal=0;\cdot)$ is not small, and destructive interference does not occur. Fortunately, this does not happen for a generic $x_0$.

The above argument applies when $x_0\in\s$ (see case (A) at the end of the introduction). If $x_0\not\in\s$, and there is a line through $x_0$ which is tangent to $\s$ (case (B) in the introduction), then the argument is somewhat similar. One of the differences between the cases is that the function $g$ is now expressed in terms of $\hal$ rather than $\tal$. Also, estimates for $g(\hal;\cdot)$ and its derivative are different from those for $g(\tal;\cdot)$. 

If lines through $x$ are never tangent to $\s$ (case (C)), no destructive interference needs to be considered. In this case, \eqref{extra terms I} follows from the smallness of $A_m(\al,\e)$. Interestingly, in all three cases, 
level sets of $H_0$ appear in the proofs in an essential way. They are used to estimate $g$ and its derivative. Nevertheless, the proof of case (C) follows a somewhat different overall logic than the proofs of cases (A) and (B) (see remark~\ref{rem orders} at the end of section~\ref{sec: sum_k II}). 

Estimates for $g(\tal_k;\cdot)$ and its derivative in case (A) are obtained in section~\ref{sec: g-ests}, and the sum with respect to $k$ is estimated in section~\ref{sec: sum_k I}. Likewise, estimates for $g(\hal_k;\cdot)$ and its derivative in case (B) are obtained in section~\ref{sec: rem prep}, and the sum with respect to $k$ is estimated in section~\ref{sec: sum_k II}. As was mentioned in the introduction, the proof of case (C) is in appendix~\ref{sec new ker prf}.

\section{Preparation}\label{sec:prep}

In sections \ref{sec:prep}--\ref{sec: sum_k I}, $x_0=y(0)$, so the function $R_1$ in \eqref{bigR1} becomes
\be\label{bigR1}
R_1(\theta,\al)=\vec\al\cdot \left(y(\theta)-y(0)\right).
\ee
The convexity of $\s$ and our convention imply that the nonzero vector $(y(\theta)-y(0))/\theta$ rotates counterclockwise as $\theta$ increases from $-a$ to $a$. Thus, for each $\theta\in[-a,a]$ there is $\al=\CA_1(\theta)\in(-\pi/2,\pi/2)$ such that $\vec\al(\CA_1(\theta))\cdot \left(y(\theta)-y(0)\right)\equiv0$, and the function $\CA_1$ is injective. By continuity, $\CA_1(0):=0$. Define
\be\label{Omega def}
\Omega:=\text{ran}\CA_1.
\ee
Clearly, $\Omega\subset(-a,a)$, and the inverse of $\CA_1(\theta)$ is smooth and well-defined on $\Omega$. In what follows we need rescaled versions of $R_1$ and $\CA_1$:
\be\label{bigR}
R(\tilde\theta,\tilde\al):=R_1(\theta,\al)/\e,\ 
\CA(\tilde\theta):=\CA_1(\theta)/\e^{1/2}.
\ee
For simplicity, the dependence of $R$ and $\CA$ on $\e$ is omitted from notation. 

\begin{definition}
We say $f(x)\asymp g(x)$ for $x\in U\subset \br^n$ if there exist $c_{1,2}>0$ such that
\be
c_1\le f(x)/g(x)\le c_2\text{ if }g(x)\not=0 \text{ and }f(x)=0\text{ if }g(x)=0
\ee 
for any $x\in U$.
\end{definition}

\begin{lemma}\label{lem:aux props}
One has 
\be\label{r props III}
\CA(\tilde\theta)\asymp\tilde\theta,\ \CA'(\tal)\asymp1,\ |\theta|\le a,|\al|\le\pi/2;\ \max_{|\tth|\le a\e^{-1/2}} |\CA(\tilde\theta)/\tilde\theta|<1,
\ee
and
\be\label{r props II}
\begin{split}
&R(\tilde\theta,\tilde \al)\asymp \tilde\theta(\CA(\tilde\theta)-\tilde\al),\ \pa_{\tilde\theta} R(\tilde\theta,\tilde \al)\asymp \tilde\theta-\tilde\al,\ \pa_{\tilde\al} R(\tilde\theta,\tilde \al)=O(\tilde\theta),\\ 
&\pa_{\tal}^2 R(\tilde\theta,\tilde \al),\, \pa_{\tilde\theta}^2 R(\tilde\theta,\tilde \al),\, \pa_{\tilde\theta}\pa_{\tilde\al} R(\tilde\theta,\tilde \al)=O(1),\quad |\theta|\le a,|\al|\le\pi/2.
\end{split}
\ee
\end{lemma}
%

Then, from \eqref{A simpl orig},
\be\label{A-simpl-ps}
\begin{split}
\e^{1/2}A_m(\al,\e)
=&\int_{\br}\text{sgn}(\hat t)\sum_n\int_{U_n}\tilde\psi_m\left(R(\tilde\theta,\tilde\al)+h(\theta,\al)\right)F(\theta,\e \hat t)\dd\tilde\theta\dd\hat t.
\end{split}
\ee
Recall that $\theta$ is a function of $\tth$ in the arguments of $h$ and $F$. Clearly,
\be\label{hF est}
h,F=O(1),\ \pa_{\tilde\theta} h,\pa_{\tilde\al} h,\pa_{\tilde\theta} F=O(\e^{1/2}),\quad |\theta|\le a,|\al|\le\pi/2,
\ee
uniformly with respect to all variables. Here $\al$ is a function of $\tal$.

Fix some small $\de>0$ and define three sets 
\be\label{three Theta sets}
\Xi_1:=[-\de,\de],\ \Xi_2:=\{\tilde\theta:\,|\CA(\tth)-\tal|\le\de\},\
\Xi_3:=[-a\e^{-1/2},a\e^{-1/2}]\setminus(\Xi_1\cup\Xi_2),
\ee
and the associated functions:
\be\label{three g}
\begin{split}
g_l(\tilde\al;m,\e,\hat t,\check x):=&\sum_n\int_{U_n\cap\Xi_l}\tilde\psi_m\left(R(\tilde\theta,\tilde\al)+h(\theta,\al)\right)F(\theta,\e \hat t)\dd\tilde\theta,\ l=1,2,3.
\end{split}
\ee
If $\al\not\in\Omega$, we assume $\Xi_2=\varnothing$ and $g_2(\tilde\al;m,\e,\hat t,\check x)=0$.
To simplify notations, the arguments $m,\e,\hat t$, and $\check x$ of $g_l$ are omitted in what follows, and we write $g_l(\tilde\al)$. In view of \eqref{W def orig}, $g=g_1+g_2+g_3$.


\section{Estimates for $g_{1,2,3}$}\label{sec: g-ests}
Using Lemma~\ref{lem:aux props} and \eqref{hF est} and estimating the model integral $\int (1+\tth^2((\tth/2)-\tal)^2)^{-1}\dd\tth$ over various domains, it is straightforward to conclude that
\be\label{g123 est}
g_{1,2}(\tilde\al)=\rho(m)O(|\tilde\al|^{-1}),\quad g_{3}(\tilde\al)=\rho(m)O(|\tilde\al|^{-2}),\quad \tal\to\infty,\ |\al|\le\pi/2.
\ee
Similarly, estimating the model integral $\int |\tth|(1+\tth^2((\tth/2)-\tal)^2)^{-1}\dd\tth$ implies
\be\label{der g13 est}
\pa_{\tilde\al}g_{1,3}(\tilde\al)=\rho(m)O(|\tilde\al|^{-1}),\quad \tal\to\infty,\ |\al|\le\pi/2.
\ee

Thus, it remains to estimate $\pa_{\tilde\al}g_2$. We assume $\al\in\Omega$, because $g_2(\tal)=0$ if $\al\not\in\Omega$. Change variables $\tilde\theta\to r=R(\tilde\theta,\tilde\al)$ in \eqref{three g}, so that $\tth=\Theta(r,\tal)$:
\be\label{g2 via r}\begin{split}
g_2(\tilde\al)=&\sum_n\int_{R_n\cap[r_{\text{mn}},r_{\text{mx}}]} \tilde\psi_m(r+h(\theta,\al))\frac{F(\theta,\e\hat t)}{|(\pa_{\tth} R)(\Theta(r,\tal),\tilde\al)|}\dd r,\\
R_n:=&R(U_n,\tilde\al),\ r_{\text{mn}}:=R(\CA^{-1}(\tal-\de),\tilde\al),\ 
r_{\text{mx}}:=R(\CA^{-1}(\tal+\de),\tilde\al),
\end{split}
\ee
where $\theta=\e^{1/2}\Theta(r,\tal)$. 
If $R$ is decreasing in $\tth$ and $r_{\text{mn}}>r_{\text{mx}}$, the domain in \eqref{g2 via r} is understood as $R_n\cap[r_{\text{mx}},r_{\text{mn}}]$. Denote also $r_n:=R(u_n,\tilde\al)$ and $v_n:=\CA(u_n)$. Clearly, $r_n$'s are the endpoints of $R_n$'s: $R_n=[r_{2n},r_{2n+1}]$ or $R_n=[r_{2n+1},r_{2n}]$ depending on whether $R(\tth,\tal)$ is increasing or decreasing as a function of $\tth$.

\begin{lemma}\label{lem: more props}
For $\al\in\Omega$, $|\tal|\ge c$, one has
\be\label{needed ests II}
\begin{split}
& r_n\asymp \tal(\CA(u_n)-\tal),\ \pa_{\tilde\al} r_n\asymp-\tal\quad\text{if}\quad |v_n-\tal|\le \de;\\
&r_{\text{mn}}\asymp -\tilde\al,\ r_{\text{mx}}\asymp \tilde\al,\ \pa_{\tal} r_{\text{mn}},\pa_{\tal} r_{\text{mx}}=O(1),\\ 
&\pa_{\tilde\theta} R(\tth,\tal)\asymp\tilde\al\quad \text{if}\quad   |\CA(\tth)-\tal|\le \de;
\end{split}
\ee
and
\be\label{ests II-1}
\pa_{\tilde\al}\Theta(r,\tal)=-\left.\frac{\pa_{\tal} R(\tth,\tal)}{\pa_{\tth} R(\tth,\tal)}\right|_{\tth=\Theta(r,\tal)}
=O(1)\quad \text{if}\quad  r\in[r_{\text{mn}}, r_{\text{mx}}].
\ee
\end{lemma}

In particular, $\pa_{\tilde\theta} R(\tth,\tal)\not=0$ if $|\CA(\tth)-\tal|\le \de$, $\al\in\Omega$, and $|\tal|\ge c$, so the change of variables in \eqref{g2 via r} is justified. Differentiating \eqref{g2 via r} and using \eqref{hF est} and Lemmas~\ref{lem:aux props}, \ref{lem: more props} gives 
\be\label{g2 final ests}\begin{split}
|\pa_{\tilde\al}g_2(\tilde\al)|\le &c\frac{\rho(m)}{|\tal|}\biggl[\sum_{n:|v_n-\tal|\le\de} \frac{|\tal|}{1+\tal^2(v_n-\tal)^2}+\frac{1}{|\tal|^2}+\e^{1/2}+\frac{1}{|\tal|}\biggr]\\
\le& c\rho(m)\left[\sum_{n:|v_n-\tal|\le\de}\frac{1}{1+\tilde\al^2(v_n-\tal)^2}+\frac{1}{\tilde\al^2}\right].
\end{split}
\ee
Here we have used that $\e^{1/2}\tal=O(1)$. 
To summarize, we have 
\be\label{g final ests}
\begin{split}
&|g(\tilde \al)|\le c\rho(m)(1+|\tilde \al|)^{-1},\ |\al|\le\pi/2;\\
&|\pa_{\tilde\al}g(\tilde \al)|\le c\rho(m)\left[\sum_{n:|v_n-\tal|\le \de}\frac{1}{1+\tilde\al^2(v_n-\tal)^2}+\frac{1}{1+\tilde\al^2}\right],\ \al\in\Omega;\\
&|\pa_{\tilde\al} g(\tilde \al)|\le c\rho(m)(1+|\tilde\al|)^{-1},\ \al\in[-\pi/2,\pi/2]\setminus\Omega;
\end{split}
\ee
where $g$ is given by \eqref{W def orig}.

\section{Summation with respect to $k$}\label{sec: sum_k I}
\subsection{Preliminary results}
The goal of this section is to estimate the sum in \eqref{W def orig}. This is done by breaking up the interval $[-\pi/2,\pi/2]$ into a union of smaller intervals and estimating individually the sums over each of these intervals. Denote
\be\label{dot pr fns}
\phi(\al)=m\kappa r_x \sin(\al-\al_x),\  \vartheta(\hal)=-m \vec\al\cdot x_0/\e,\
\mu=-\kappa r_x \sin\al_x,
\ee
where $x_0=r_x\vec\al(\al_x)$. Recall that \eqref{scales} is always assumed. Clearly, $\vartheta'(\hal)\equiv\phi(\al)$, $\phi(0)=m\mu$, and $\mu=-\kappa\theta_0^\perp\cdot x_0$. Without loss of generality we may assume $m\ge 1$. We begin by estimating the top sum in \eqref{W def orig}. There are four cases to consider: $m\ge 1$ or $m\le -1$ combined with $\al_k\in[0,\pi/2]$ or $\al_k\in[-\pi/2,0]$. We will consider only one case:  $m\ge 1$ and $\al_k\in[0,\pi/2]$, the other three cases are completely analogous.

Let $\al_*>0$ be the smallest angle such that $\phi'(\al_*)=0$, i.e. $\al_*=\al_x+(\pi/2)$ (mod $\pi$). 
Assumption P2 implies $\al_*\not=0$. Otherwise, $\vec\theta_0\cdot x_0=0$, $\vec\theta_0\cdot y'(0)=0$, and $x_0=y(0)$ imply that the line through the origin and $x_0$ is tangent to $\s$.
By assumption P4, $\al_*\not=\pi/2$. If not, $\kappa\theta_0^\perp\cdot x_0=0$ is rational. 

\begin{figure}[h]
{\centerline
{\epsfig{file=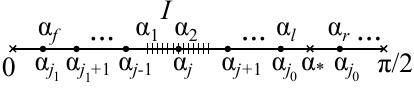, width=9cm}}}
\caption{Illustration of the interval $[0,\pi/2]$ with various angles used in the derivation of the estimates.}
\label{fig:alphas}
\end{figure}

Let $\al_{s,m}$ satisfy  
\be\label{alsm}
\phi(\al_{s,m})=s,\ |s|\le m\kappa r_x,\ s\in(1/2)\BZ,\ \al_{s,m}\in [0,\pi/2].
\ee
See Figure~\ref{fig:alphas}, where $\al_{s,m}$ are shown as thick dots for integer values of $s$ (and without the subscript $m$).
If $\al_*\in(0,\pi/2)$, then for some $s$ there may be two solutions: $\al_{s,m}\in (0,\al_*)$ and $\al_{s,m}\in (\al_*,\pi/2]$. If $\al_*\not\in(0,\pi/2)$, there is at most one solution for each $s$. Note that the solution(s) may exist only for some of the indicated $s$. 
By assumption P3, $\kappa r_x$ is irrational, so $\al_{s,m}\not=\al_*$ for any $s,m$. By assumption P4, $\phi(0)$ is irrational, so $\al_{s,m}\not=0$ for any $s,m$. 

Consider an interval $I\subset[0,\pi/2]$ and its rescaled versions
\be\label{rescaled ints}
\hat I:=(1/\Delta\al)I,\ \tilde I:=\e^{-1/2}I.
\ee
In view of \eqref{extra terms I} and \eqref{W def orig}, consider the expression
\be\label{W def I}
W_m(I):=\e^{1/2} \left|\sum_{\tal_k\in \tilde I}e\left(-m q_k\right) \left[g(\tilde\al_k)e\left(-m \vec\al_k\cdot \check x\right)\right]\right|,\ q_k:=\frac{\vec\al_k\cdot x_0}\e.
\ee
The dependence of $W_m(I)$ on $\e$ and $\hat t$ is omitted from notations. The goal is to estimate the sum $\sum_{m\not=0}W_m([0,\pi/2])$. Obviously, $W_m([0,\pi/2])\le\sum_j W_m(I_j)$ if $\cup_j I_j=[0,\pi/2]$.
From \eqref{g final ests}, 
\be\label{simple W est}
W_m([0,\pi/2])\le O(\e^{1/2})\rho(m)\sum_{k=0}^{1/\e}(1+\e^{1/2}k)^{-1}
=\rho(m)O(\ln(1/\e)),
\ee
because there are $O(1/\e)$ values $\al_k\in[0,\pi/2]$. Therefore,
\be\label{sum W est}
\sum_{|m|\ge c\e^{-1/2}}W_m([0,\pi/2])\le O(\ln(1/\e)) \sum_{m\ge c\e^{-1/2}}\rho(m)
=O(\e\ln(1/\e))
\ee
if $\rho(m)=O(|m|^{-3})$ and $c>0$. Thus, in what follows, we will assume $\e^{1/2}|m|\le c$ for some small $c>0$. 

Next we investigate the individual $W_m(I)$ for smaller intervals $I$. For this we need a partial integration identity \cite[p. 89]{Huxley1996} (written in a slightly different form):
\be\label{pii}
\sum_{k=K_1}^{K_2} G(k)\Phi(k)=G(K_2)\sum_{k=K_1}^{K_2} \Phi(k)-\int_{K_1}^{K_2}G'(\tau)\sum_{k=K_1}^{\tau} \Phi(k)\dd \tau,
\ee
where $G(\tau)$ is continuously differentiable on the interval $[K_1,K_2]$. Here and throughout the paper, $\sum_{k=c_1}^{c_2}$, where $c_{1,2}$ are not necessarily integers, denotes the sum over $k\in[c_1,c_2]$. Hence
\be\label{pii-ineq}
\left|\sum_{k=K_1}^{K_2} G(k)\Phi(k)\right|\le \left(|G(K_2)|+\int_{K_1}^{K_2}|G'(\tau)|\dd \tau\right)\max_{k'\in[K_1,K_2]}\left|\sum_{k=K_1}^{k'} \Phi(k)\right|.
\ee
The following result is also needed \cite[p. 7]{Graham1991}:
\begin{theorem}[Kusmin-Landau inequality] If $\vartheta(\tau)$ is continuously differentiable, $\vartheta'(\tau)$ is monotonic, and $\langle \vartheta'(\tau)\rangle\ge \la>0$ on an interval $\hat I$, then there exists $c>0$ (independent of $\vartheta$ and $\hat I$) such that 
\be\label{KL ineq}
\left|\sum_{n\in \hat I}e(\vartheta(n))\right|\leq c/\la.
\ee
\end{theorem}

In view of \eqref{W def I}, \eqref{pii}, define
\be\label{Gh}
\begin{split}
G(\hal):=g(\tal)e\left(-m \vec\al \cdot \check x\right),\, \Phi(\hal):=e(\vartheta(\hal)),
\end{split}
\ee
where $\vartheta(\hal)$ is defined in \eqref{dot pr fns}.

\begin{lemma}\label{sum int bnd} One has for any $b,L\ge 0$, $L=O(1)$, $[b,b+L]\subset\e^{-1/2}\Omega$:
\be\label{strange int stp1}
\int_b^{b+L}
\sum_{n:|v_n-\tal|\le \de}\frac{1}{1+\tilde\al^2(v_n-\tilde\al)^2}\dd \tilde\al
= O((1+b)^{-1}).
\ee
\end{lemma}

\begin{corollary}\label{cor: int bnd} One has for any $b,L\ge 0$, $L=O(1)$, $b+L\le\e^{-1/2}\pi/2$:
\be\label{full gder int}
\int_{b}^{b+L} |\pa_{\tilde\al} g(\tilde \al)|\dd \tal
= \rho(m)O((1+b)^{-1}).
\ee
\end{corollary}
The assertion is obvious because (a) the second term in brackets in \eqref{g final ests} is $O(\tal^{-2})$, and (b) by \eqref{g final ests}, $\pa_{\tilde\al} g(\tilde \al)$ satisfies the same estimate as in \eqref{strange int stp1} if $\tal\not\in \e^{-1/2}\Omega$ (with $\rho(m)$ accounted for).

Pick any interval $I\subset[0,\pi/2]$ such that $|\tilde I|=O(1)$. From \eqref{Gh}, 
\be\label{var bnd stp1}
\int_{\hat I}|G'(\hal)|\dd \hal
\le c\int_{\tilde I}\left(|\pa_{\tilde\al}g(\tilde\al)|+\e^{1/2}|m||g(\tilde\al)|\right)\dd \tilde\al.
\ee
Also, by \eqref{g final ests} and \eqref{full gder int}:
\be\label{var bnd stp2}
\begin{split}
&|G(\hal_0)|=|g(\tal_0)|=\rho(m)O((1+\tal_0)^{-1}),\\
& \int_{\tilde I}|g(\tal)|\dd \tal,
\int_{\tilde I}|\pa_{\tilde\al}g(\tilde\al)|\dd \tilde\al
= \rho(m)O((1+\tal_0)^{-1})\text{ for any } \al_0\in I.
\end{split}
\ee
Hence, 
\be\label{G bnd comb}
\begin{split}
\max_{\hal\in\hat I}|G(\hal)|+\int_{\hat I}|G'(\hal)|\dd \hal
=\rho(m)O((1+\tal_0)^{-1})\text{ for any } \al_0\in I,
\end{split}
\ee
where we have used that $\e^{1/2}m\le c$. 

We suppose that $\al=\al_*$ is a local maximum of $\phi(\al)$. The case when $\al_*$ is a local minimum is analogous. 

Let $[0,\al_f]$ be the shortest interval such that $\phi(\al_f)$ is an integer, i.e. $\phi(\al_f)\in\BZ$ and $\phi(\al)\not\in\BZ$ for any $\al\in[0,\al_f)$. Similarly, let $[\al_l,\al_r]$ be the shortest interval centered at $\al_*$ such that $\phi(\al_l)=\phi(\al_r)\in\BZ$ and $\phi(\al)\not\in\BZ$ for any $\al\in(\al_l,\al_r)$. Clearly, $\phi(\al_l)=\phi(\al_r)=\lfloor\phi(\al_*)\rfloor$. Set
\be\label{I0 def}
\Iz:=[0,\pi/2]\setminus ([0,\al_f]\cup[\al_l,\al_r]).
\ee
Thus, $W_m([0,\pi/2])\le W_m(\Iz)+W_m([0,\al_f])+W_m([\al_l,\al_r])$.

\subsection{Towards estimation of $W_m(\Iz)$}\label{ssec:bulk}
In this subsection we assume $\al_f<\al_l$, i.e. $\lceil\phi(0)\rceil\le\lfloor\phi(\al_*)\rfloor$, because otherwise $\Iz=\varnothing$, and $W_m(\Iz)=0$.

Pick any interval $I:=[\al_1,\al_2]\subset \Iz$ (if one exists) such that $2\phi(\al_{1,2})\in\BZ$, $|\phi(\al_2)-\phi(\al_1)|=1/2$, and $2\phi(\al)\not\in\BZ$ for any $\al\in(\al_1,\al_2)$. By construction, $\phi(\al)$ is monotone on $I$. Let $j$ be the integer value in the pair $\phi(\al_1),\phi(\al_2)$, i.e. $\al_1=\al_{j,m}$ or $\al_2=\al_{j,m}$ (see \eqref{alsm}). The other value in the pair is $j-1/2$ or $j+1/2$. 

Generically, one has
\be\label{cos equiv}
|\cos(\al_1-\al_x)|\asymp |\cos(\al_2-\al_x)|.
\ee
The only exceptions are the two cases when $[\al_1,\al_2]$ is close to $\al_*$: $\al_2=\al_l$ and $\al_1=\al_r$. In these cases one of the expressions in \eqref{cos equiv} can be arbitrarily close to zero (as $m\to\infty$), while the other can stay away from zero. If $\al_2=\al_l$, then $\phi(\al_2)=j$ and $\phi(\al_1)=j-1/2$. If $\al_1=\al_r$, then $\phi(\al_1)=j$ and $\phi(\al_2)=j-1/2$. Therefore, when $m\gg1$,
\be\label{cos equiv special}\begin{split}
&|\cos(\al_1-\al_x)|> |\cos(\al_{j,m}-\al_x)| \text{ if } \al_{j,m}=\al_l,\\
&|\cos(\al_2-\al_x)|> |\cos(\al_{j,m}-\al_x)| \text{ if } \al_{j,m}=\al_r.
\end{split}
\ee
Clearly, $\al_2-\al_1=O(m^{-1/2})$. Away from a neighborhood of $\al_*$, this difference is actually $O(1/m)$. To cover all the cases we use a more conservative estimate. 

Subdivide $\tilde I$ into $N$ subintervals of length $\asymp1$ (see Figure~\ref{fig:alphas}): 
\be\label{Ijl sub}\begin{split}
&\tilde I=\cup_{n=0}^{N-1} [\tal_{j,m}-(n+1)L,\tal_{j,m}-nL],\ \tal_{j,m}-NL=\tal_1 \text{ if }\al_2=\al_{j,m},\\
&\tilde I=\cup_{n=0}^{N-1} [\tal_{j,m}+nL,\tal_{j,m}+(n+1)L],\ \tal_{j,m}+NL=\tal_2 \text{ if }\al_1=\al_{j,m}.
\end{split}
\ee
Clearly, $\tal_2-\tal_1\ge c/(\e^{-1/2}m)$, so the requirement $\e^{-1/2}m\le c$ implies $\tal_2-\tal_1\ge c$. Hence we can choose $N=1$ and $L=\tal_2-\tal_1$ if $\tal_2-\tal_1<1$, and $N=\lfloor\tal_2-\tal_1\rfloor$ and $L=(\tal_2-\tal_1)/N$ if $\tal_2-\tal_1\ge 1$.

One has:
\be\label{sin ineq}\begin{split}
|\sin\al-\sin\al_1|\ge (\al-\al_1)\min(|\cos\al_1|,|\cos\al_2|),\ \al\in [\al_1,\al_2],\\
|\sin\al_2-\sin\al|\ge (\al_2-\al)\min(|\cos\al_1|,|\cos\al_2|),\ \al\in [\al_1,\al_2],
\end{split}
\ee
for any $\al_1<\al_2$ such that $\sin\al\not=0$ on the interval $(\al_1,\al_2)$. The statement is immediate in view of the mean value theorem and the monotonicity of $\cos\al$ on the interval $[\al_1,\al_2]$.\

By construction, $j$ is the integer closest to $\vartheta'(\hal)$ if $\hal\in\hat I$:
\be\label{dist to int}
\langle \vartheta'(\hal)\rangle = |\vartheta'(\hal)-j|=|\vartheta'(\hal)-\vartheta'(\hal_{j,m})|\text{ if }
\hal\in \hat I.
\ee
From \eqref{cos equiv}, \eqref{cos equiv special}, \eqref{sin ineq}, and \eqref{dist to int}
\be\begin{split}
&\langle \vartheta'(\hal)\rangle\ge c(\e^{1/2}nL) m\kappa r_x|\cos(\al_{j,m}-\al_x)|=c\e^{1/2}n((m\kappa r_x)^2-j^2)^{1/2}\\ 
&\text{if }\tal\in[\tal_{j,m}-(n+1)L,\tal_{j,m}-nL]\text{ or }\tal\in[\tal_{j,m}+nL,\tal_{j,m}+(n+1)L].
\end{split}
\ee
This follows from the top inequality in \eqref{sin ineq} if $\al_1=\al_{j,m}$, and from the bottom one -- if $\al_2=\al_{j,m}$. Also,
\be\label{al ineq}
\tal\ge c\tal_{j,m}\text{ if }\al\in I.
\ee
Indeed, notice that $\tal\ge \tal_{j,m}$ if $\al_1=\al_{j,m}$. If $\al_2=\al_{j,m}$, suffices it to assume that $m\gg 1$ is large enough. From $\al_*\not=0$, $\al_*\not\in[0,c']$ for some $c'>0$. By construction, $|\phi(\al_1)-\phi(0)|>0.5$, $|\phi(\al_2)-\phi(0)|>1$, and $|\phi(\al_2)-\phi(\al_1)|=1/2$. If $m\gg1$, \eqref{al ineq} is obvious if $\al_2>c'$, and it follows from the mean value theorem if $[\al_1,\al_2]\subset [0,c']$.
%

The partial integration identity \eqref{pii} and the Kusmin-Landau inequality \eqref{KL ineq} imply
\be\label{Wjm bulk est}\begin{split}
W_m(I)
&\le c\e^{1/2}\rho(m)\left(\frac{\e^{-1/2}}{\tal_{j,m}}+\frac{\e^{-1/2}}{\tal_{j,m}((m\kappa r_x)^2-j^2)^{1/2}}\sum_{n=1}^N\frac1n\right)\\
&\le  c\frac{\rho(m)}{\tal_{j,m}}\left(1+\frac{\ln(1/(\e m))}{((m\kappa r_x)^2-j^2)^{1/2}}\right).
\end{split}
\ee
The first term in parentheses on the first line in \eqref{Wjm bulk est} bounds the contribution from the subinterval $[\tal_{j,m}-L,\tal_{j,m}]$ or $[\tal_{j,m},\tal_{j,m}+L]$ (depending on the case), which is adjacent to $\tal_{j,m}$. Since $\phi(\al_{j,m})=j$, we cannot use the Kusmin-Landau inequality, so $W_m(\cdot)$ for this subinterval is estimated directly from \eqref{W def I} using the top line in \eqref{g final ests} and \eqref{al ineq}.

Clearly, $\Iz$ can be represented as a union of intervals $I=[\al_1,\al_2]$ of the kind considered in this subsection. Summing the estimates in \eqref{Wjm bulk est} for all $I\subset\Iz$ to obtain a bound for $W_m(\Iz)$ is done in section~\ref{ssec:comb} below. Therefore, it is left to consider $W_m([\al_l,\al_r])$ and $W_m([0,\al_f])$. 

\subsection{Estimation of $W_m([\al_l,\al_r])$}\label{ssec: al*}
Suppose first that $\al_*\in(0,\pi/2)$. 
Since $\al_*$ is a local maximum, $\phi(\al)$ is increasing on $[\al_l,\al_*]$ and decreasing  - on $[\al_*,\al_r]$. When $m\gg1$ is sufficiently large, we have $0<\al_l<\al_*<\al_r\le\pi/2$. 

Suppose $\al_l>0$ (see appendix~\ref{sec: spec cases} if this does not hold). Then $\tal\asymp\e^{-1/2}$ if $\tal\in[\tal_l,\tal_r]$. Split $[\tal_l,\tal_*]$ into $N=\lfloor\tal_*-\tal_l\rfloor$ subintervals of length $L\asymp 1$:
\be\label{I*l sub v0}\begin{split}
&[\tal_l,\tal_*]=\cup_{n=0}^{N-1} [\tal_l+nL,\tal_l+(n+1)L],\ \tal_l+NL=\tal_*.
\end{split}
\ee
Since $\tal_*-\tal_l \asymp [\{m\kappa r_x\}/(\e m)]^{1/2}$, we have
\be\label{bigN*}
c[\langle m\kappa r_x\rangle/(\e m)]^{1/2}\le \tal_*-\tal_l \le c/(\e m)^{1/2}.
\ee

Applying \eqref{G bnd comb} to each of the subintervals in \eqref{I*l sub v0} gives the estimate $\rho(m)O((1+\tal_*)^{-1})=\rho(m)O(\e^{1/2})$.
Also,
\be\label{thpr*}\begin{split}
\langle \vartheta'(\hal)\rangle&\ge \min\left(\frac{nL}{\tal_*-\tal_l}\{\phi(\al_*)\},
1-\{\phi(\al_*)\}\right)\\
&\ge \min\left(cn(\e m \langle m\kappa r_x\rangle)^{1/2},
\langle m\kappa r_x\rangle\right),\ 
\tal\in[\tal_l+nL,\tal_l+(n+1)L].
\end{split}
\ee
Completely analogous estimates hold for $[\tal_*,\tal_r]$ if $\al_r\le\pi/2$ (see also appendix~\ref{sec: spec cases}). Therefore,
\be\label{Wjm est stp1}\begin{split}
W_m([\al_l,\al_r])\le &c\e^{1/2}\rho(m)\left(1+\frac{\e^{1/2}}{(\e m)^{1/2}\langle m\kappa r_x\rangle}+\frac1{(m\langle m\kappa r_x\rangle)^{1/2}}\sum_{n=1}^{N-1}\frac1n\right)\\
\le & c\e^{1/2}\rho(m)\left(\frac1{m^{1/2}\langle m\kappa r_x\rangle}+\frac{\ln\left(1/(\e m)\right)}{(m\langle m\kappa r_x\rangle)^{1/2}}\right),\ \tal_*-\tal_l\ge 1.
\end{split}
\ee
The first term in parentheses on the first line in \eqref{Wjm est stp1} corresponds to the subinterval $[\tal_l,\tal_l+L]$, since $\phi(\tal_l)\in \BZ$ and its contribution is estimated directly from \eqref{W def I} using the top line in \eqref{g final ests}.

If $\tal_*-\tal_l < 1$, we can estimate $W_m([\al_l,\al_r])$ directly from \eqref{W def I}. There are $O(\e^{-1/2})$ terms in the sum, each of them is $O(\tal_*^{-1})=O(\e^{1/2})$, so 
\be\label{Wjm est v2}
W_m([\al_l,\al_r])\le c\e^{1/2}\rho(m),\ \tal_*-\tal_l < 1. 
\ee

\subsection{Estimation of $W_m([0,\al_f])$}\label{ssec:al0}
Since $\al_*$ is a local maximum, $\phi(\al)$ is increasing on $[0,\al_*]$. 
If there is no $\al\in (0,\al_*)$ such that $\phi(\al)\in\BZ$, then $\al_l< 0$, and this case is addressed in appendix~\ref{sec: spec cases}. Therefore, in this subsection we assume that $\phi(\al_{j_1,m})=j_1$, where $j_1:=\lceil \phi(0) \rceil$, for some $\al_{j_1,m}\in (0,\al_*)$. Clearly, $\al_f=\al_{j_1,m}$.

Split $[0,\tal_{j_1,m}]$ into $N=\lfloor\tal_{j_1,m}\rfloor$
intervals of length $L\asymp 1$:
\be\label{I0 sub}\begin{split}
[0,\tal_{j_1,m}]=\cup_{n=0}^{N-1} [nL,(n+1)L],\ NL=\tal_{j_1,m}, \text{ if }\tal_{j_1,m}\ge 2.
\end{split}
\ee
Since $\tal_{j_1,m}\asymp (1-\{m\mu\})/(\e^{1/2} m)$, we have
\be\label{tal bnds}
c\langle m\mu\rangle/(\e^{1/2} m)\le\tal_{j_1,m}\le c/(\e^{1/2} m).
\ee
Applying \eqref{G bnd comb} to the interval $\hat I_n=(\kappa\e^{1/2})^{-1}[nL,(n+1)L]$ gives
\be\label{G bnd comb 0}
\max_{\hal\in\hat I_n}|G(\hal)|+\int_{\hat I_n}|G'(\hal)|\dd \hal
=\rho(m)O((1+n)^{-1}).
\ee
By the bottom line in \eqref{sin ineq}, 
\be
j_1-\phi(\al)\ge \e^{1/2}(N-(n+1))L(m\kappa r_x)\min(|\cos(-\al_x)|,|\cos(\al_{j_1,m}-\al_x)|),
\ee
if $\tal\in [nL,(n+1)L]$, therefore
\be\label{deriv est 1}
\langle \vartheta'(\hal)\rangle\ge \min\left(\{ m\mu\},c\e^{1/2}m(N-(n+1))\right),\ \hal\in\hat I_n,
\ee
because $|\cos(-\al_x)|,|\cos(\al_{j_1,m}-\al_x)|\asymp1$. Combining the inequalities, using \eqref{tal bnds}, and simplifying gives
\be\label{Wj1m est}\begin{split}
W_m([0,\al_{j_1,m}])\le &c\e^{1/2}\rho(m)\biggl(\frac{\e^{-1/2}}{\tal_{j_1,m}}+\frac{1}{\{m\mu\}}\sum_{n=0}^{N-2}\frac1{1+n}\\
&\qquad
+\frac1{\e^{1/2}m}\sum_{n=0}^{N-2}\frac1{(1+n)(N-(n+1))}\biggr)\\
\le & c\e^{1/2}\rho(m)\frac{m+\ln\left(1/(\e^{1/2} m)\right)}{\langle m\mu\rangle},\ \tal_{j_1,m}\ge 2.
\end{split}
\ee
The first term in parentheses on the first line in \eqref{Wj1m est} bounds the contribution from the last subinterval $\tilde I_{N-1}=[\tal_{j_1,m}-L,\tal_{j_1,m}]$. Since $\phi(\al_{j_1,m})=j_1$, we cannot use the Kusmin-Landau inequality, so $W_m(I_{N-1})$ is estimated directly from \eqref{W def I} using the top line in \eqref{g final ests}. From these two equations we get also
\be\label{Wj1m est v2}
W_m([0,\al_{j_1,m}])\le c\rho(m),\ \tal_{j_1,m}<2.
\ee


\subsection{Combining all the estimates.}\label{ssec:comb}
We begin by summing \eqref{Wjm bulk est} over all the intervals $I=[\al_1,\al_2]\subset\Iz$ in order to finish estimating $W_m(\Iz)$. The analysis in section~\ref{ssec:bulk} shows that $W_m(I)$ admits the same bound \eqref{Wjm bulk est} regardless of whether $\al_1=\al_{j,m}$ or  $\al_2=\al_{j,m}$. Hence we need to sum the right-hand side of \eqref{Wjm bulk est} over all integers $j\in\phi(\Iz)$. Recall that $\al_{j,m}$ denote the angles such that $\phi(\al_{j,m})=j$ (see \eqref{alsm}). We need to distinguish two cases: $0<\al_{j,m}\le\min(\al_l,\pi/2)$ and $\al_r\le\al_{j,m}<\pi/2$. The latter case may occur only if $\al_*<\pi/2$.

First, suppose $0<\al_{j,m}\le\min(\al_l,\pi/2)$. Denote $j_1:=\lceil\phi(0)\rceil$, $j_0:=\lfloor\phi(\al_*)\rfloor$. Then
\be
\al_{j,m}\ge [\langle m\mu\rangle +(j-j_1)]/(m\kappa r_x),\ j_1\le j \le j_0,
\ee
and
\be\label{sum Wjm est}\begin{split}
W_m([\al_f,\al_l])
\le c\e^{1/2}\rho(m)\biggl[&\frac{m}{\langle m\mu\rangle}\left(1+\frac{\ln(1/(\e m))}m\right)
+\left(1+\frac{\ln(1/(\e m))}{m^{1/2}\langle m\kappa r_x\rangle^{1/2}}\right)\\
&+\sum_{j=j_1+1}^{j_0-1}\frac{m}{j-j_1}\left(1+\frac{\ln(1/(\e m))}{m^{1/2}(m\kappa r_x-j)^{1/2}}\right)\biggr].
\end{split}
\ee
The first term in brackets corresponds to $j=j_1$, the second term -- to $j=j_0$, and the sum -- to all $j_1<j<j_0$. If $j_0-j_1=1$, the sum is assumed to be zero.

If $\al_l>\pi/2$, $W_m([\al_f,\pi/2])$ is assumed instead of $W_m([\al_f,\al_l])$.
Clearly, $m\kappa r_x-j\ge 1$ if $j_1<j<j_0$.
Simplifying and keeping only the dominant terms we have
\be\label{sum Wjm final}\begin{split}
W_m([\al_f,\al_l])
&\le c\e^{1/2}\rho(m)\left(\frac{m+\ln(1/(\e m))}{\langle m\mu\rangle}+\frac{\ln(1/(\e m))}{m^{1/2}\langle m\kappa r_x\rangle^{1/2}}\right).
\end{split}
\ee


If $\al_r<\pi/2$ and $\al_r<\al_{j,m}\le\pi/2$, then $\tal_{j,m}\asymp \e^{-1/2}$, so \eqref{Wjm bulk est} gives
\be\label{sum Wjm est v2}\begin{split}
W_m([\al_r,\pi/2])
\le & c\e^{1/2}\rho(m)\left(m+\frac{\ln(1/(\e m))}{m^{1/2}}\left[\langle m\kappa r_x\rangle^{-1/2}+m^{1/2}\right] \right).
\end{split}
\ee
Here we used that there are $O(m)$ distinct values of $j$. Clearly, the estimate in \eqref{sum Wjm final} dominates the one in \eqref{sum Wjm est v2}. This implies that $W_m(\Iz)$ satisfies the estimate in \eqref{sum Wjm final}. 

We have
\be\label{two rates}
\rho(m)=O(m^{-\bt}),\ \langle m\mu\rangle,\langle m|x_0|\rangle=O(m^{-\eta}),\ m\to\infty,
\ee
for any $\eta>\eta_0$ (see conditions P3, P4). By summing each of the estimates \eqref{Wjm est stp1} (contribution of $[\al_l,\al_r]$), \eqref{Wj1m est} (contribution of $[0,\al_f]$), and \eqref{sum Wjm final} (contribution of $\Iz$) with respect to $m$ from 1 to $\infty$, we see that the dominating term is $\e^{1/2}\sum_{m=1}^\infty\rho(m) m/\langle m\mu\rangle$ (cf. \eqref{Wj1m est} and \eqref{sum Wjm final}). The series converges if $\bt>\eta_0+2$. 

Next consider special cases. Comparing \eqref{Wjm est v2} and \eqref{Wj1m est v2}, we see that the latter grows faster as $m\to\infty$. From \eqref{tal bnds} and \eqref{Wj1m est v2},
\be
\sum_{m\ge 1,\tal_{j_1,m}<2}W_m([0,\al_{j_1,m}])\le c
\sum_{\substack{m\ge 1\\ \langle m\mu\rangle/(\e^{1/2} m)\le c}}\rho(m)=O(\e^{(\bt-1)/(2(\eta+1))})=O(\e^{1/2})
\ee
if $\bt\ge \eta_0+2$. 

The contribution of the exceptional cases that take place for finitely many $m$ (see appendix~\ref{sec: spec cases}) is of order $O(\e^{1/2}\ln(1/\e))$. Hence we proved that 
\be\label{sum Wjm bnd}
\sum_{1\le |m|\le O(\e^{-1/2})}W_m([-\pi/2,\pi/2])=O(\e^{1/2}\ln(1/\e))\text{ if }\bt>\eta_0+2.
\ee

Given that $\hat t$ is confined to a bounded interval and \eqref{sum Wjm bnd} is uniform with respect to $\hat t$ (see the paragraph following \eqref{W def orig}), we prove \eqref{extra terms I} in the case $x_0\in\s$.

\section{Remote singularity: Preparation}\label{sec: rem prep}

\subsection{Preliminaries} 
Now we consider the case $x_0\not\in\s$. The following result is proven in \cite{Katsevich2021b}.

\begin{lemma}\label{lem:partial res} Pick $x_0\not\in\s$ such that no line through $x_0$, which intersects $\s$, is tangent to $\s$. This includes the endpoints of $\s$, in which case the one-sided tangents to $\s$ are considered. Under the assumptions of Theorem~\ref{main-res}, one has
\be\label{extra terms top}
\frec(x)=O(\e^{1/2}\ln(1/\e)),\ \e\to0,
\ee
uniformly with respect to $x$ in a sufficiently small neighborhood of $x_0$.
\end{lemma}

To make this paper self-contained, a slightly simplified proof of the lemma is in Appendix~\ref{sec new ker prf}. By the above lemma and a partition of unity only a small neighborhood of a point of tangency can be considered. Therefore, in this section also we can assume that $\s$ is as short ($a>0$ as small) as we like. 

In this section we use $R_1$ in its original form (see \eqref{bigR1 orig})
\be\label{*bigR1}
R_1(\theta,\al)=\vec\al\cdot \left(y(\theta)-x_0\right),
\ee
because $x_0\not=y(0)$.
Let $\al=\CA_1(\theta)$: $[-a,a]\to[-\pi/2,\pi/2]$ be the function such that $\vec\al(\CA_1(\theta))\cdot(y(\theta)-x_0)\equiv0$. Suppose, for example, that $x_0$ is on the side of $\s$ for which $\CA_1(\theta)\ge0$. The case when $\CA_1(\theta)\le0$ is analogous.
In contrast with Section~\ref{sec:prep}, $\CA_1$ is now quadratic near $\theta=0$. In what follows we need rescaled versions of the functions $R_1$ and $A_1$:
\be\label{*bigRA}
R(\tilde\theta,\hat\al):=R_1(\theta,\al)/\e,\ 
\CA(\tilde\theta):=\CA_1(\theta)/\Delta\al.
\ee
As usual, the dependence of $R$ and $\CA$ on $\e$ is omitted from notation for simplicity.   
\begin{lemma}\label{lem:aux props*}
One has 
\be\label{*r props I}
\CA(\tilde\theta)\asymp\tilde\theta^2,\ \pa_{\tilde\theta} \CA(\tilde\theta)\asymp\tilde\theta,\ |\theta|\le a;
\ee
and
\be\label{*r props II}
\begin{split}
&R(\tilde\theta,\hat \al)\asymp \CA(\tilde\theta)-\hat\al,\ \pa_{\tilde\theta} R(\tilde\theta,\hat \al)\asymp \tilde\theta-\kappa\e^{1/2}\hat\al,\ \pa_{\hat\al} R(\tilde\theta,\hat \al)=O(1),\\ 
&\pa_{\tilde\theta}^2 R(\tilde\theta,\hat \al)=O(1),\ \pa_{\tth}\pa_{\hal} R(\tilde\theta,\hat \al)=O(\e^{1/2}),\
|\theta|\le a,|\al|\le\pi/2.
\end{split}
\ee
\end{lemma}
The proof of this and all other lemmas in this section are in Appendix~\ref{sec: *lems prep}.
From \eqref{A-simpl},
\be\label{*A-simpl-ps}
\begin{split}
\e^{1/2}A_m(\al,\e)=&\int_{-a\e^{-1/2}}^{a\e^{-1/2}}\int_{\br}\tilde\psi_m\left(R(\tilde\theta,\hat\al)+h(\theta,\al)\right)
F(\theta,\e \hat t)\chi(\hat t,H_0(\tilde\theta))\dd\hat t\dd\tilde\theta\\
=&\int_{\br}\text{sgn}(\hat t)\sum_n\int_{U_n}\tilde\psi_m\left(R(\tilde\theta,\hat\al)+h(\theta,\al)\right)F(\theta,\e \hat t)\dd\tilde\theta\dd\hat t.
\end{split}
\ee
Clearly,
\be\label{*hF est}
h,F=O(1),\ \pa_{\tilde\theta} h,\pa_{\tilde\theta} F=O(\e^{1/2}),\ \pa_{\hat\al} h=O(\e),\quad |\theta|\le a,|\al|\le\pi/2,
\ee
uniformly with respect to all variables. By \eqref{*A-simpl-ps},
\be\label{* g def}
\begin{split}
&A_m(\al,\e)=\e^{-1/2}\int_{\br}\text{sgn}(\hat t) g(\hat\al;m,\e,\hat t,\check x)\dd\hat t,\\
&g(\hat\al;m,\e,\hat t,\check x):=\sum_n\int_{U_n}\tilde\psi_m\left(R(\tilde\theta,\hat\al)+h(\theta,\al)\right)F(\theta,\e \hat t)\dd\tilde\theta.
\end{split}
\ee
As usual, the arguments $m,\e,\hat t$, and $\check x$ of $g$ are omitted in what follows, and we write $g(\hat\al)$. 

Our goal is to estimate the sum in \eqref{extra terms I}. To do that we first simplify the sum by reducing the range of indices $k$ and simplifying the expression for $A_m$. Note that we no longer assume $m\ge 1$.

\subsection{Simplification of the sum \eqref{extra terms I}}
From \eqref{*r props I}--\eqref{*hF est} and \eqref{four-coef-bnd}, it is easy to obtain $g(\hat\al)=\rho(m)O(|\hat\al|^{-3/2})$, $\hat\al\to-\infty$. This implies
\be\label{neg-k}
\begin{split}
\Delta\al\sum_m \sum_{\hat\al_k\le \hat \al_*} |A_m(\al_k,\e)|=O(\e^{1/2})
\end{split}
\ee
for any fixed $\hat \al_*>0$. The meaning of $\al_*$ here is different from that in Section~\ref{sec: sum_k I}. Similarly to \eqref{Omega def}, introduce the set $\Omega=\text{ran}\CA_1$ (with $\hat\Omega=(1/\Delta\al) \Omega$ according to our usual convention).
We will show that the sum over $\al_k\in [0,\pi/2]\setminus\Omega$ makes only a negligible contribution to $\frec$. By \eqref{neg-k}, the contribution of negative $\al_k$ and any finite number of $\al_k>0$ can be ignored. Introduce the variable $\tau$:
\newcommand{\Amx}{\CA_{\text{mx}}}
\be\label{tau def}
\hat\al=\Amx+\tau,\ \Amx:=\max(\CA(-a\e^{-1/2}),\CA(a\e^{-1/2})).
\ee
Clearly, $\hat\Omega=\text{ran}\CA=[0,\Amx]$. 

\begin{lemma}\label{lem:large pos k} One has
\be\label{large pos k}
\begin{split}
\Delta\al\sum_m \sum_{\hat\al_k > \Amx} |A_m(\al_k,\e)|
=O(\e\ln(1/\e)).
\end{split}
\ee
\end{lemma}

\begin{lemma}\label{lem:g1 bnd} One has
\be\label{g1 bound}
g_1(\hat\al):=\sum_n\int_{\substack{\tth\in U_n\\ |R(\tilde\theta,\hat\al)|\ge \de \hat\al^{1/2}}}\tilde\psi_m\left(R(\tilde\theta,\hat\al)+h(\theta,\al)\right)F(\theta,\e \hat t)\dd\tilde\theta=\rho(m)O(\hat\al^{-1})
\ee
as $\hat\al\to+\infty$ for any $\de>0$ as small as we like. 
\end{lemma}

Equation \eqref{g1 bound} implies
\be\label{pos-k 1}
\Delta\al\sum_m \sum_{\hat\al_k\in[\hat\al_*,\Amx]} |A_m^{(1)}(\al_k,\e)|
=O(\e^{1/2}\ln(1/\e)),
\ee
where $A_m^{(1)}$ is obtained by the top line in \eqref{* g def} with $g$ replaced by $g_1$.

%

The only remaining contribution to $\frec$ comes from 
\be\label{g2 def}
g_2(\hat\al):=\sum_n\int_{\substack{\tth\in U_n \\ |R(\tilde\theta,\hat\al)|< \de \hat\al^{1/2}}}\tilde\psi_m\left(R(\tilde\theta,\hat\al)+h(\theta,\al)\right)F(\theta,\e \hat t)\dd\tilde\theta,\ \hat\al\in\hat\Omega.
\ee
A simple calculation shows that $g_2(\hat\al)=O(\hat\al^{-1/2})$, $\hat\al\to\infty$. By \eqref{*r props I} and \eqref{*r props II},
\be\label{del theta}\begin{split}
&\tth-\CA^{-1}(\hat\al)=O(\hat\al^{1/2}/\CA'(\CA^{-1}(\hat\al)))
=O(1),\quad\theta-\Theta(\hat\al)=O(\e^{1/2}),\\
&\Theta(\hat\al):=\e^{1/2}\CA^{-1}(\hat\al),\ \hat\al\to\infty,\ \hat\al\in\hat\Omega,\
|R(\tilde\theta,\hat\al)|< \de \hat\al^{1/2}.
\end{split}
\ee
Even though $\CA^{-1}$ is two-valued, \eqref{del theta} holds regardless of whether $\theta>0$ and $\text{ran} \CA^{-1}\subset [0,\infty)$ or $\theta<0$ and $\text{ran} \CA^{-1}\subset (-\infty,0]$. Hence, by \eqref{four-coef-bnd} we can replace $F(\theta,\e \hat t)$ and $h(\theta,\al)$ with $F(\Theta(\hat\al),\e \hat t)$ and $h(\Theta(\hat\al),\al)$, respectively, in \eqref{g2 def}:
\be\label{g2 approx}
\begin{split}
g_2(\hat\al)=&g_2^+(\hat\al)+g_2^-(\hat\al)+\rho(m)O((\e/\hat\al)^{1/2}),\ \hat\al\in \hat\Omega,\\
g_2^\pm(\hat\al)=&F(\Theta(\hat\al),\e \hat t)\sum_n\int_{\substack{\pm\tth>0,\tth\in U_n, \\ |R(\tilde\theta,\hat\al)|< \de \hat\al^{1/2}}}\tilde\psi_m\left(R(\tilde\theta,\hat\al)+h(\Theta(\hat\al),\al)\right)\dd\tilde\theta.
\end{split}
\ee
The superscript $'+'$ is taken if $\text{ran}\CA^{-1}=[0,\infty)$, and $'-'$ -- if $\text{ran}\CA^{-1}=(-\infty,0]$. 
The same convention is assumed in what follows. 
In particular, the domain of integration in \eqref{g2 approx} is a subset of $(0,\infty)$ when $g_2^+$ is computed, and a subset of $(-\infty,0)$ -- when $g_2^-$ is computed.
Omitting the big-$O$ term in \eqref{g2 approx} leads to:
\be\label{pos-k 2}
\Delta\al\sum_m \sum_{\hat\al_k\in[\hat\al_*,\Amx]} |A_m^{(1)}(\al_k,\e)-(A_m^{(2+)}(\al_k,\e)+A_m^{(2-)}(\al_k,\e))|=O(\e^{1/2}),
\ee
where $A_m^{(2\pm)}$ are obtained by the top line in \eqref{* g def} with $g$ replaced by $g_2^{\pm}$, respectively. Due to this simplification, we can consider
\be\label{g3 def}
g_3^{\pm}(\hat\al):=F\sum_n\int_{\substack{\pm\tth>0,\tth\in U_n, \\ |R(\tilde\theta,\hat\al)|< \de \hat\al^{1/2}}}\tilde\psi_m\left(R(\tilde\theta,\hat\al)+h\right)\dd\tilde\theta,
\ee
where $F=F(\Theta(\hat\al),\e \hat t)$, $h=h(\Theta(\hat\al),\al)$ are uniformly bounded and independent of $\tth$. As was done before, set $r_n(\hat\al):=R(u_n,\hat\al)$. By shifting the index of $u_n$ if necessary, we may suppose that $u_0$ is the smallest nonnegative $u_n$. This means that $u_n\ge 0$ if $n\ge 0$ and $u_n<0$ if $n<0$. In this case, $u_n\asymp n$ if $|u_n|\ge c$ by assumption $H3$. Changing variables $\tilde\theta\to r=R(\tilde\theta,\hat\al)$ gives:
\be\label{g3pl r}
g_3^\pm(\hat\al):=F\sum_{\pm n > 0}\int_{\substack{r\in R_n \\ |r|\le\de\hat\al^{1/2}}} \tilde\psi_m(r+h)\frac{\dd r}{|(\pa_{\tth} R)(\tth,\hat\al)|},
\ee
where 
$\tth$ is a function of $r$ and $\hat\al$. As usual, $+n>0$ in $g_4^+$, and $-n>0$ -- in $g_4^-$.
The change of variables is justified, because $(\pa_{\tth} R)(\tth,\hat\al)\not=0$ on the integration domain. Indeed, $\tilde\theta$ is bounded away from zero on the domain, and the following result holds. 

\begin{lemma}\label{lem: nozero} One has
\be\label{tth equiv}
(\pa_{\tth} R)(\tth,\hat\al)
\asymp\tth \quad\text{if}\quad |R(\tilde\theta,\hat\al)|< \de\hat\al^{1/2}, \hat\al\in[\hat\al_*,\Amx].
\ee
\end{lemma}


Further simplification is achieved by replacing $\tth$ with $\CA^{-1}(\hat\al)$ in the argument of $\pa_{\tth} R$. From \eqref{*r props II} and \eqref{del theta}
\be\label{g4pl r}\begin{split}
g_3^\pm(\hat\al)=&g_4^\pm(\hat\al)+\rho(m)O(\hat\al^{-1}),\\
g_4^\pm(\hat\al):=&\frac{F}{|(\pa_{\tth} R)(\CA^{-1}(\hat\al),\hat\al)|}\sum_{\pm n>0}\int_{\substack{r\in R_n\\|r|\le\de\hat\al^{1/2}}} \tilde\psi_m(r+h)\dd r,\\ 
F=&F(\Theta(\hat\al),\e \hat t),\ h=h(\Theta(\hat\al),\e\hat\al).
\end{split}
\ee
Neglecting the big-$O$ term in \eqref{g4pl r} leads to a term of magnitude $O(\e^{1/2}\ln(1/\e))$ in $\frec$. 

Here is a summary of what we obtained so far:
\be\label{*summary}
\begin{split}
&\Delta\al\sum_m \sum_{\hat\al_k\not\in[\hat\al_*,\Amx]} |A_m(\al_k,\e)|=O(\e^{1/2}),\\
&\Delta\al\sum_m \sum_{\hat\al_k\in[\hat\al_*,\Amx]} |A_m(\al_k,\e)-(A_m^{+}(\al_k,\e)+A_m^{-}(\al_k,\e))|=O(\e^{1/2}\ln(1/\e)),
\end{split}
\ee
where (cf. \eqref{* g def})
\be\label{last A}
A_m^{\pm}(\al,\e):=\e^{-1/2}\int_{\br}g_4^{\pm}(\hat\al;m,\e,\hat t,\check x)\dd\hat t.
\ee
The first line in \eqref{*summary} follows from \eqref{neg-k} and \eqref{large pos k}. The second line follows from \eqref{pos-k 1}, \eqref{pos-k 2}, and the comment following \eqref{g4pl r}.

Finally, we need the following result.
\begin{lemma}\label{g4 bnds} One has
\begin{align}\label{g2 asymp}
&|g_4^\pm(\hat\al)|\le c \rho(m) \left(\hat\al^{-1/2}\left(1+r_{\text{min}}(\hal)\right)^{-1}+\hat\al^{-1}\right),\\
\label{g2derasymp}
&|\pa_{\hat\al}g_4^\pm(\hat\al)|\le c \rho(m) \left(\frac1{\hat\al^{1/2}(1+r_{\text{min}}^2(\hal))}
+\frac{\e^{1/2}}{\hal}\right),
\end{align}
where $\hat\al\in[\hat\al_*,\Amx]$, and
\be\label{rmin def}
r_{\text{min}}(\hal):=\min_{\pm n>0}|\CA(u_n)-\hat\al|.
\ee
\end{lemma}
To clarify, in the estimate for $g_4^+$, the minimum in \eqref{rmin def} is over $n>0$, and in the estimate for $g_4^-$ -- over $n<0$. 

\section{Summation with respect to $k$}\label{sec: sum_k II}

Denote $v_n:=\CA(u_n)$ and consider the intervals
\be\label{rn ints}\begin{split}
&V_n:=[v_{n-0.5},v_{n+0.5}],n>0,\ V_n:=[v_{n+0.5},v_{n-0.5}],n<0,\\ 
&v_{n+0.5}:=(v_n+v_{n+1})/2,\ n\in\BZ.
\end{split}
\ee
Since the function $\CA(\tth)$ in sections~\ref{sec: rem prep} and \ref{sec: sum_k II} is different from the one in sections~\ref{sec:prep}--\ref{sec: sum_k I}, the $v_n$'s in \eqref{rn ints} are different from the $v_n$'s in sections~\ref{sec:prep}--\ref{sec: sum_k I}.

Clearly, the estimates in Lemma~\ref{g4 bnds} increase if extra points are (formally) added to the list of $u_n$'s. If $u_n\not\asymp n$ (i.e., there are too few $u_n$'s), we can always add more points to the $u_n$'s and enumerate them so that the enlarged collection satisfies $u_n\asymp n$. This property is assumed in what follows.

\begin{lemma}\label{lem: Vn props}
If $V_n\subset[\hat\al_*,\Amx]$, one has
\be\label{vVn props}
v_n\asymp n^2,\ |V_n|=|v_{n+0.5}-v_{n-0.5}|\asymp |n|;\
\hal\asymp v_n,\ r_{min}(\hal)=|v_n-\hal|\text{ if }\hal\in V_n
\ee
\end{lemma}

The proof of the lemma is immediate using Lemma~\ref{lem:aux props*}, \eqref{rn ints}, and that $u_n\asymp n$. 

\begin{lemma}\label{lem: on Vn}
For all $V_n$ such that $V_n\subset[\hat\al_*,\Amx]$, one has as $n\to\infty$:
\be\label{g4 bnd endpt}
|g_4^\pm(v_{n+1/2})|=\rho(m)O(n^{-2}),
\ee
and
\be\label{g4 ints}
\int_{V_n}|g_4^\pm(\hat\al)|\dd \hat\al = \rho(m)O(|n|^{-1}\ln |n|),\
\int_{V_n}|\pa_{\hat\al}g_4^\pm(\hat\al)|\dd \hat\al =\rho(m)O(|n|^{-1}).
\ee
\end{lemma}

Arguing similarly to \eqref{simple W est}, \eqref{sum W est}, eq. \eqref{g4 ints} implies 
\be\label{*simple W est}
\begin{split}
\e^{1/2}\sum_{\hat\al_k\in[\hat\al_*,\Amx]} |A_m^{\pm}(\al_k,\e)|
&\le \rho(m)\sum_{n=1}^{O(\e^{-1/2})}n^{-1}\ln n
=\rho(m)O(\ln^2(1/\e)).
\end{split}
\ee
Here we used Lemma~\ref{lem: Vn props} and the following two arguments. (i) Since $\hat\al_*>0$ can be taken as large as we want, we can select $\hal_*=\hal_*(\hat t,\e)$ so that (a) $c_1\le \hat\al_*\le c_2$ for some fixed $c_{1,2}>0$ and all $\hat t$ and $\e>0$, and (b) $\hal_*=v_{n+0.5}$ for some $n$. (ii) If $v_{n_l}$ is the last of the $v_n$'s in the interval $[\hat\al_*,\Amx]$, the sum over $\hal_k\in [v_{n_l},\Amx]$ can be estimated directly. By \eqref{vVn props}, $\Amx-\hal_*=O(\Amx^{1/2})$, so the number of the additional $\hal_k$'s is $O(\Amx^{1/2})$. By \eqref{g2 asymp}, $g_4^\pm(\hal)=\rho(m)O(\Amx^{-1/2})$, $\hal_k\in [v_{n_l},\Amx]$. Therefore the contribution of $[v_{n_l},\Amx]$ is $O(\rho(m))$, which is absorbed by the right-hand side of \eqref{*simple W est}. 

From \eqref{*simple W est},
\be\label{*sum W est}
\begin{split}
&\Delta\al\sum_{|m|\ge\ln(1/\e)} \sum_{\hat\al_k\in[\hat\al_*,\Amx]} |A_m^\pm(\al_k,\e)|
=O(\e^{1/2})
\end{split}
\ee
if $\rho(m)=O(|m|^{-3})$. Thus, in what follows, we will assume $|m|\le\ln(1/\e)$.

Let  $\al_m^{(0)}$ be the smallest positive root of the equation $|\phi(\al)-\phi(0)|=\langle \phi(0)\rangle/2$. 
Recall that the function $\phi$ depends on $m$. By \eqref{two rates} and the restriction on $m$ (recall that $\phi(0)=m\mu$),
\be\label{hal lower bnd}
\langle m\mu\rangle \ge c (\ln(1/\e))^{-\eta},\quad
\hal_m^{(0)}\ge c(\e\ln^{\eta}(1/\e))^{-1},
\ee
for any $\eta>\eta_0$. Find $n$ such that $\hal_m^{(0)}\in V_{n}$ and set $\hal_m:=v_{n+0.5}$. By construction and \eqref{vVn props}, 
\be
0<\hal_m-\hal_m^{(0)}=O(v_{n}^{1/2})=O((\hal_m^{(0)})^{1/2})=O(\e^{-1/2}),
\ee
therefore $\al_m-\al_m^{(0)}=O(\e^{1/2})$.

In this section we use only two intervals 
\be\label{two ints}
\hat I_1:=[\hat\al_*,\hal_m],\
\hat I_2:=[\hal_m,\Amx],\ m\not=0.
\ee 
Since $\hal_*=O(1)$, due to \eqref{hal lower bnd} we can assume $\hal_m\ge \hat\al_*$. The derivations for $g_4^\pm$ are the same, so we drop the superscript. Similarly to \eqref{W def I}, define
\be\label{*W def}
\begin{split}
W_m(I):=&\e^{1/2} \biggl|\sum_{\hat\al_k\in \hat I}e\left(-m q_k\right) \left[g_4(\hat\al_k)e\left(-m \vec\al_k\cdot \check x\right)\right]\biggr|.
\end{split}
\ee
The dependence of $W_m(I)$ on $\e$ and $\hat t$ is omitted from notations.

Following the method in Section~\ref{sec: sum_k I}, we use the partial integration identity \eqref{pii-ineq} and the Kusmin-Landau inequality \eqref{KL ineq}. The definitions of $\Phi$ and $\vartheta$ are the same as before, and the definition of $G$ (cf. \eqref{Gh}) is modified slightly 
\be\label{*Gh}
\begin{split}
G(\hal):=g_4(\hal)e\left(-m \vec\al\cdot \check x\right).
\end{split}
\ee

We begin by applying \eqref{pii-ineq} to the first interval, so we select $[K_1,K_2]$ as the interval such that $k\in [K_1,K_2]$ is equivalent to $\hal_k\in\hat I_1$.
By \eqref{hal lower bnd} and the choice of $\al_m^{(0)}$, $\al_m$,
\be\label{*spacing 0I}
\hal_m\ge c\langle m\mu\rangle/(\e|m|),\quad 
\langle \vartheta'(\hal)\rangle\ge c\langle m\mu\rangle,\, \hal\in\hat I_1,\ m\not=0.
\ee
From \eqref{*Gh}, 
\be\label{*var bnd stp1}
\begin{split}
\int_{\hat I_1}|G'(\hal)|\dd \hal
\le c\int_{\hat I_1}\left(|\pa_{\hat\al}g_4(\hat\al)|+\e|m||g_4(\hat\al)|\right)\dd \hat\al.
\end{split}
\ee

Our construction ensures that $\hat I_1=\cup V_n$, where the union is taken over all $V_n\subset \hat I_1$. 
There are $N=O((\hal_m^{(0)})^{1/2})$ intervals $V_n$ such that $V_n\subset\hat I_1$. By \eqref{g4 bnd endpt}, \eqref{g4 ints},
\be\label{*var bnd stp2}
\begin{split}
&|G(K_2)|\le c |g_4(v_{N+0.5})| \le c \rho(m) N^{-2},\\
&\int_{\hat I_1}|g_4(\hat\al)|\dd \hat\al
\le c\rho(m)\sum_{n=1}^{N}O(n^{-1}\ln n)=\rho(m) O(\ln^2 N),\\
&\int_{\hat I_1}|\pa_{\hat\al}g_4(\hat\al)|\dd \hat\al
\le c\rho(m)\sum_{n=1}^{N}n^{-1}=\rho(m) O(\ln N).
\end{split}
\ee
Hence, 
\be\label{*G bnd comb}
\begin{split}
|G(K_2)|+\int_{\hat I_1}|G'(\hal)|\dd \hal
&=\rho(m) O(\ln N).
\end{split}
\ee
Here we have used that $|m|\le\ln(1/\e)$ implies $\e|m|\ln N<1$ for $\e>0$ sufficiently small. For the same reason we can assume $N\ge 1$. By \eqref{*spacing 0I}, \eqref{*G bnd comb}, \eqref{pii-ineq}, and the Kusmin-Landau inequality,
\be\label{*W1m est}
W_m(I_1)\le c\frac{\e^{1/2}\rho(m)}{\langle m\mu\rangle}\ln\left(\frac{\langle m\mu\rangle}{\e|m|}\right),\ m\not=0.
\ee
The sum over the remaining $V_n\subset \hat I_2$, $N\le n \le O(\e^{-1/2})$, can be estimated easily without utilizing exponential sums. By construction, the left endpoint of $\hat I_2$ coincides with $v_{n+0.5}$ for some $n$. As was established following \eqref{*simple W est}, the contribution of $\hal_k$ beyond the last $V_n\subset \hat I_2$ is $O(\rho(m))$.
Therefore Lemma~\ref{lem: on Vn} implies 
\be\label{I1l bnd}\begin{split}
\sum_{\hat\al_k\in\hat I_2}|g_4(\hat\al_k)|&\le 
c\int_{\hat I_2}|g_4(\hat\al)|\dd \hat\al \le c\sum_{V_n\subset\hat I_2}\int_{V_n}|g_4(\hat\al)|\dd \hat\al +O(\rho(m))\\
&\le c\rho(m)\sum_{n=N}^{O(\e^{-1/2})}n^{-1}\ln n \le c\rho(m)(\ln^2(\e^{-1/2})-\ln^2(N)) \\
&=\rho(m) O\left(\ln(1/\e)\ln(|m|/\langle m\mu\rangle)\right).
\end{split}
\ee
With some abuse of notation, in the two integrals on the top line above, we integrate the upper bound for $g_4$ obtained in \eqref{g2 asymp} (with $\hal$ replaced by $v_n$ by Lemma~\ref{lem: Vn props}). This bound has better monotonicity properties (i.e., it can be made monotone within $V_n$ on each side of $v_n$). Otherwise, we would not be able to estimate the sum in terms of an integral. We also used that $\hal_{k+1}-\hal_k=1$. 

Thus,
\be\label{W2m est}
W_m(I_2)=\rho(m)\ln(|m|/\langle m\mu\rangle) O(\e^{1/2}\ln(1/\e)),\ m\not=0.
\ee
Comparing \eqref{*W1m est} and \eqref{W2m est} with \eqref{sum Wjm final} we see that the case $x_0\not\in\s$ adds no additional restrictions on $\rho(m)$. Similarly to the end of section~\ref{sec: sum_k I}, we use here that the estimates \eqref{*W1m est} and \eqref{W2m est} are uniform with respect to $\hat t$ and $\hat t$ is confined to a bounded interval.

\begin{remark} \label{rem orders}
\rm{Note that the order of operations in the proof of Lemma~\ref{lem:partial res} (see appendix~\ref{sec new ker prf}) is as follows:
\be
\sum_{\al_k}(\cdot)\to \int (\cdot)\dd\hat t\to \sum_m \int(\cdot) \dd \tth.
\ee
In section~\ref{sec:beg proof}--\ref{sec: sum_k II} the order is different:
\be
\int (\cdot)\dd\tth \to \sum_{\al_k}(\cdot)\to \sum_m \int(\cdot) \dd \hat t.
\ee}
\end{remark}

\section{The case $m=0$}\label{m=0 case}
Here we prove \eqref{extra terms II}. Begin with the case $x_0=y(0)$. By \eqref{A simpl orig}, \eqref{W def orig},
\be\label{2nd term stp1}
\begin{split}
&J_\e:=\sum_{|\al_k|\le \pi/2}
\int_{\al_k-\Delta\al/2}^{\al_k+\Delta\al/2}  \left|A_0(\al,\e)-A_0(\al_k,\e)\right| \dd\al\le O(\e^{1/2})\int J_\e(\hat t)\dd \hat t, \\
&J_\e(\hat t):=\sum_{|\al_k|\le \pi/2}\max_{|\tal-\tal_k|\le \Delta\tal/2}  \left|\pa_{\tal}g(\tal)\right| \Delta\tal,\ 
\Delta\tal:=\kappa\e^{1/2}.
\end{split}
\ee
As is seen from \eqref{g final ests}, the only term that requires careful estimation is given by
\be\label{Je1 st1}
\begin{split}
J_\e^{(1)}(\hat t):=&\Delta\tal\sum_{\al_k\in\Omega}\max_{|\tal-\tal_k|\le \Delta\tal/2}  \sum_{n:|v_n-\tal|\le \de}g_1(\tal,v_n),\\ 
g_1(\tal,v_n):=&\left[1+v_n^2(v_n-\tal)^2\right]^{-1}.
\end{split}
\ee
Here we used that $1+\tilde\al^2(v_n-\tal)^2\asymp 1+v_n^2(v_n-\tal)^2$ if $|v_n-\tal|\le \de$.
Replacing the inner sum with a larger sum over $n:|v_n-\tal_k|\le 2\de$ (using that $\de+(\Delta\al/2)\le 2\de$)
we can write
\be\label{Je1 st2}
\begin{split}
J_\e^{(1)}(\hat t)\le \Delta\tal\sum_{|n|\le O(\e^{-1/2})}\sum_{k:|v_n-\tal_k|\le 2\de} \max_{|\tal-\tal_k|\le \Delta\tal/2}  g_1(\tal,v_n).
\end{split}
\ee
Since $\pa_{\tal}g_1(\tal,v_n)<0$, $\tal>v_n$, 
\be\label{Je1 st3}
\begin{split}
&\sum_{k:v_n\le \tal_k\le v_n+2\de} \max_{|\tal-\tal_k|\le \Delta\tal/2}  g_1(\tal,v_n)
\le \sum_{0\le k\Delta\tal \le 2\de}   g_1(v_n+k\Delta\tal,v_n)\\
&\le 1+\frac1{\Delta\tal}\int_{v_n}^{\infty}g_1(\tal,v_n)\dd\tal\le c\left(1+\frac1{\Delta\tal(1+|v_n|)}\right).
\end{split}
\ee
The same argument applies to the left of $v_n$, so by \eqref{vVn props}
\be\label{Je1 st4}
\begin{split}
J_\e^{(1)}(\hat t)\le&\sum_{|n|\le O(\e^{-1/2})}\left(O(\e^{1/2})+(1+n^2)^{-1}\right)=O(1).
\end{split}
\ee

Consider now the remaining terms in \eqref{g final ests}. Define similarly to \eqref{Je1 st1}:
\be\label{Je23}
\begin{split}
J_\e^{(2)}(\hat t):=&\Delta\tal\sum_{\al_k\in\Omega}\max_{|\tal-\tal_k|\le \Delta\tal/2} (1+\tal^2)^{-1},\\ 
J_\e^{(3)}(\hat t):=&\Delta\tal\sum_{\al_k\in[-\pi/2,\pi/2]\setminus\Omega}\max_{|\tal-\tal_k|\le \Delta\tal/2} (1+|\tal|)^{-1}.
\end{split}
\ee
Arguing analogously to \eqref{Je1 st3}, \eqref{Je1 st4} we get $J_\e^{(l)}(\hat t)=O(1)$, $l=2,3$. When estimating $J_\e^{(3)}$ we use that the summation is over $\al_k$ which satisfy $c\le |\al_k|\le \pi/2$.
Combining with \eqref{Je1 st4} and substituting into \eqref{Je1 st1} leads to $J_\e=O(\e^{1/2})$. 

Suppose now $x_0\not\in\s$. It is obvious that the continuous analogue of \eqref{*summary} works for $m=0$: \be\label{*summary cont}
\begin{split}
&\int_{\hat\al\not\in [\hat\al_*,\Amx]} |A_0(\al,\e)|\dd\al=O(\e^{1/2}),\\
&\int_{\hat\al\in [\hat\al_*,\Amx]} |A_0(\al,\e)-(A_0^{+}(\al,\e)+A_0^{-}(\al,\e))|\dd\al=O(\e^{1/2}\ln(1/\e)).
\end{split}
\ee
Hence it remains to estimate
\be\label{m=0 est}
J_\e^\pm:=\sum_{\hat\al_k\in [\hat\al_*,\Amx]}
\int_{\al_k-\Delta\al/2}^{\al_k+\Delta\al/2} |A_0^\pm(\al,\e)-A_0^\pm(\al_k,\e)|\dd\al.
\ee
By \eqref{last A},
\be\label{I est}
J_\e^\pm\le c\e^{1/2}\int J_\e^\pm(\hat t) \dd \hat t,\quad
J_\e^\pm(\hat t):=\sum_{V_n \subset [\hat\al_*,\Amx]} \sum_{\hat\al_k\in V_n} 
\max_{|\hal-\hal_k|\le 1/2}  \left|\pa_{\hal}g_4^\pm(\hal)\right|+O(1).
\ee
Here $O(1)$ is the contribution of $\hal_k$ beyond the last $V_n\subset[\hat\al_*,\Amx]$ (see the argument following \eqref{*simple W est}). Clearly, we can assume $\hal_*>1/2$. By \eqref{g2derasymp}, \eqref{rmin def}, and \eqref{vVn props}, it is easy to see that
\be\label{I est part}\begin{split}
&\sum_{\hat\al_k\in V_n} 
\max_{|\hal-\hal_k|\le 1/2}  \left|\pa_{\hal}g_4^\pm(\hal)\right|\le
c\sum_{\hat\al_k\in V_n} 
\left(\frac1{v_n^{1/2}(1+(v_n-\hal_k)^2)}
+\frac{\e^{1/2}}{v_n}\right).
\end{split}
\ee
Arguing similarly to \eqref{Je1 st3} we conclude that the left-hand side of \eqref{I est part} is bounded by the same estimate as the integral of $\left|\pa_{\hal}g_4^\pm(\hal)\right|$ in \eqref{g4 ints}.
Combining \eqref{m=0 est}--\eqref{I est part} gives the desired result:
\be\label{I est final}
J_\e^\pm\le c\e^{1/2}\sum_{n=1}^{O(\e^{-1/2})}n^{-1}=O(\e^{1/2}\ln(1/\e)).
\ee

%
%


\appendix

\section{Proofs of lemmas in sections~\ref{sec:beg proof}-\ref{sec: sum_k I}}\label{sec:loc sing proofs}

\subsection{Proof of Lemma~\ref{lem:psi psider}}
By assumption AF1, $\tilde w(\la)=O(|\la|^{-(\lceil \bt\rceil+1)})$, $\la\to\infty$, where  $\tilde w(\la)$ is the Fourier transform of $w$. If $t$ is restricted to any compact set, the estimate for $\tilde\psi_m$ holds because $\widetilde{(\CH\ik')}(\la),\tilde w(\la)=O(\rho(\la))$ (cf. assumption IK1) implies 
\be\label{convol decay}\begin{split}
\left|\int |\mu|\tilde\ik(\mu)\tilde w(\mu-\la)e^{i(\mu-\la)t}\dd\mu\right|
\le \int |\mu\tilde\ik(\mu)\tilde w(\mu-\la)|\dd\mu=O(\rho(\la)), \la=2\pi m\to\infty.
\end{split}
\ee
If $|t|\ge c$ for some $c\gg1$ sufficiently large, integrate by parts $\lceil \bt\rceil$ times and use that 
\be\label{deriv bnd}
\max_{q} |(\pa/\pa q)^{\lceil \bt\rceil}((\CH\ik')(q)w(-q-t))|=O(t^{-2}),\ t\to\infty.
\ee
The argument works, because $(\CH\ik')(q)$ is smooth in a neighborhood of any $q$ such that $w(-q-t)\not=0$. 

The estimate for $\tilde\psi_m'$ follows by differentiating \eqref{four-ser} and applying the above argument with $w$ replaced by $w'$. The argument still works because $w'\in C_0^{\lceil \bt\rceil}(\br)$ (by AF1).

\subsection{Proof of Lemma~\ref{lem:aux props}}
To prove \eqref{r props III} we write
\be\label{A stp1}
\vec\al(\CA_1(\theta))\cdot\frac{y(\theta)-y(0)}\theta\equiv
\vec\al(\CA_1(\theta))\cdot\left(y'(0)+y''(0)(\theta/2)+O(\theta^2)\right)\equiv0.
\ee
Differentiating with respect to $\theta$ and setting $\theta=0$ gives
\be\label{A stp2}
\vec\theta_0^\perp\cdot y'(0)\CA_1'(0)+\vec\theta_0\cdot y''(0)(1/2)=0.
\ee
Hence $\CA_1'(0)=1/2$ (because $\vec\al^\perp\cdot y'(\al)+\vec\al\cdot y''(\al)\equiv0$), and the desired properties of $\CA$ follow by rescaling $\theta\to\tilde\theta$ and $\CA_1\to\CA$. 

By the choice of coordinates,
\be
\vec\al\cdot \left(y(\theta)-y(0)\right)/(\theta(\CA(\theta)-\al))
\ee 
is a smooth positive function of $(\al,\theta)\in [-\pi/2,\pi/2]\times[-a,a]$. The positivity follows from the statements: (i) $\s$ is convex, (ii) by construction, $\vec\al^\perp(0)\cdot y'(0)<0$ and $\vec\al(0)\cdot y''(0)>0$, and (iii) $\vec\al\cdot(y(\theta)-y(0))$ has first order zero at $\theta=0$ if $\al\not=0$ and at $\CA_1(\theta)=\al$ if $\al,\theta\not=0$. Rescaling $\theta\to\tilde\theta$ and $\al\to\tilde\al$ proves the first property in \eqref{r props II}. 

The rest of \eqref{r props II} follows by differentiating $\vec\al\cdot \left(y(\theta)-y(0)\right)$, using that $\vec\al\cdot y'(\al)\equiv0$, and rescaling.

\subsection{Proof of Lemma~\ref{lem: more props}}
Suppose $|\tal|\ge c$. Assuming $|\CA_1(\theta)-\al|\le \de\e^{1/2}$, where $\de>0$ is sufficiently small, the properties (see \eqref{r props III})
\be
\CA_1(\theta)\asymp\theta,\quad \max_{|\theta|\le a}|\CA_1(\theta)/\theta|<1, 
\ee
imply $\theta\asymp\al$, $\theta/\al\ge c'$ for some $c'>1$, and $\pa_{\tilde\theta} R(\tth,\tal)\asymp\tilde\al$ (see \eqref{r props II}). Also, differentiating $R_1$ in \eqref{bigR1} we get $\pa_\al R_1(\theta,\al)\asymp-\theta$ and, hence, $\pa_{\tal} R(\tth,\tal)\asymp-\tal$. The properties of $r_n$ now follow immediately by rescaling and setting $\tth=u_n$. The magnitudes of $r_{\text{mn}}$ and $r_{\text{mx}}$ follow as well, because $\CA(\tth)-\tal= \pm\de$ for the corresponding $\tth$.

Denote $B(\tal):=\CA^{-1}(\tal)$. To prove the statement about $\pa_{\tal} r_{\text{mn}},\pa_{\tal} r_{\text{mx}}$ we need to show that $\pa_{\tal}R(B(\tal\pm\de),\tal)=O(1)$. Using that (i) $\pa_{\tal}\pa_{\tth} R,\pa_{\tal}^2 R$, and $B'$ are all $O(1)$ (see \eqref{r props II}, $B'=O(1)$ follows from $A'=O(1)$), and (ii) $R(B(\tal\pm\de),\tal\pm\de)\equiv0$, it is easy to see that the desired assertion holds. 

To prove \eqref{ests II-1} we differentiate $R(\Theta(r,\tal),\tal)\equiv0$ and use that $\pa_{\tth} R(\tth,\tal)\asymp\tilde\al$ and $\pa_{\tal} R(\tth,\tal)\asymp-\tilde\al$.

\subsection{Proof of Lemma~\ref{sum int bnd}}

Any $v_n$ such that $|v_n-\tal|\le\de$ for some $\tal\in[b,b+L]$ satisfies $b-\de\le v_n\le b+L+\de$. 
By Lemma~\ref{lem:aux props}, $u_n$ and $v_n$ satisfy qualitatively the same assumptions (see assumption $H3$). 
Since $L=O(1)$, there are finitely many such $v_n$. Also, $1+\tal^2(v_n-\tal)^2\asymp 1+b^2(v_n-\tal)^2$. Let $A$ denote the expression on the left side of \eqref{strange int stp1}. Then
\be\label{asymp denom}
A\le \sum_{b-\de\le v_n\le b+L+\de}\int_{\br}\frac{\dd\tal}{1+b^2(v_n-\tal)^2}=O(1/b),\ b\to\infty.
\ee

\section{Proofs of lemmas in Sections~\ref{sec: rem prep} and \ref{sec: sum_k II}}\label{sec: *lems prep}

\subsection{Proof of Lemma~\ref{lem:aux props*}}
Using that $\vec\theta_0\cdot y'(0)=0$ and $\vec\theta_0\cdot (y(0)-x_0)=0$, it is easy to see that the right-hand side of the following identity
\be\label{sin A1}
\sin(\CA_1(\theta))=\vec\theta_0\cdot \left(y(\theta)-x_0\right)/|y(\theta)-x_0|
\ee
and its first derivative are zero at $\theta=0$. Also, its second derivative at $\theta=0$ equals $\vec\theta_0\cdot y''(0)/|y(0)-x_0|$. By our choice of coordinates, this expression is positive. The properties in \eqref{*r props I} now follow from the properties of $\sin^{-1}(t)$ and by rescaling. 

To prove the first property in \eqref{*r props II}, consider the function $R_1(\theta,\al)/(\CA_1(\theta)-\al)$ and note that $\vec\al(\CA_1(\theta))$ is the only unit vector with $|\al|<\pi/2$ orthogonal to $y(\theta)-x_0$. 
Recall that $a$ is sufficiently short, so $\CA_1([-a,a])\subset[-\pi/2,\pi/2]$.
This function is clearly smooth on $[-a,a]\times[-\pi/2,\pi/2]$. The ratio is positive, because (i) $\s$ is convex, (ii) $\vec\theta_0\cdot y''(0)>0$ (cf. the proof of Lemma~\ref{lem:aux props}) and $\CA_1(\theta)\ge0$ by the assumption about $x_0$ (see the paragraph following \eqref{*bigR1}), and (iii) $R_1(\theta,\al)$ has a root of first order at $\al=\CA_1(\theta)$, $\theta\not=0$. Rescaling $\theta\to\tilde\theta$ and $\al\to\hat\al$ we finish the proof. 

Differentiating $R_1$ and using that $\vec\al\cdot y'(\al)\equiv0$ we find 
\be
\pa_\theta R_1(\theta,\al)\asymp \theta-\al;\quad \pa_\theta^2 R_1(\theta,\al),\pa_\theta\pa_\al R_1(\theta,\al),\pa_\al R_1(\theta,\al)=O(1).
\ee 
Rescaling the variables proves the rest of \eqref{*r props II}.

\subsection{Proof of Lemma~\ref{lem:large pos k}}
By \eqref{*r props II},
\be\label{large al stp1}
|g(\Amx+\tau)|\le c\rho(m)\int_{-a \e^{-1/2}}^{a \e^{-1/2}}\frac{1}{1+(\tau+(\Amx-\CA(\tth)))^2}\dd\tilde\theta,\ \tau>0.
\ee
By construction, $\Amx-\CA(\tth)\ge0$. By \eqref{*r props I}, we can change variables $r=\CA(\tth)$ (separately on $(-a\e^{-1/2},0]$ and $[0,a\e^{-1/2})$). Then $\dd\tth/\dd r\asymp \pm r^{-1/2}$ and
\be\label{large al stp2}
|g(\Amx+\tau)|\le c\rho(m)\int_{0}^{\Amx}\frac{1}{(1+\tau+\Amx-r)^2r^{1/2}}\dd r,\ \tau>0.
\ee
Here we used that $1+x^2\asymp(1+x)^2$, $x\ge 0$. Since $\Amx\asymp 1/\e$, we obtain 
\be\label{large al stp4}
\int_{0}^{\Amx/2}(\cdot)\dd r\le c \Amx^{-2}\int_{0}^{\Amx/2}\frac{\dd r}{r^{1/2}}
=O(\e^{3/2}),\ \tau>0,
\ee
and
\be\label{large al stp4}
\int_{\Amx/2}^{\Amx}(\cdot)\dd r\le c \Amx^{-1/2}\int_{0}^{\Amx/2}\frac{\dd r}{(1+\tau+r)^2}
\le c \frac{\e^{1/2}}{1+\tau},\ \tau>0.
\ee
Consequently, the sum in \eqref{large pos k} is bounded by
\be\label{large al stp6}
c\e^{1/2} \sum_{0\le k\le O(1/\e)}\left[\e^{3/2}+\frac{\e^{1/2}}{1+k}\right]=O(\e\ln(1/\e)).
\ee

\subsection{Proof of Lemma \ref{lem:g1 bnd}}
By Lemma~\ref{lem:aux props*}, there exists $\de'>0$ such that for any $0<\al\le\pi/2$, one has
\be\label{two sets}
\{|\theta|\le a: |R(\tth,\hal)|\ge\de\hal^{1/2}\}\subset \{|\theta|\le a: |\CA(\tth)-\hal|\ge\de'\hal^{1/2}\}.
\ee
This implies
\be\label{g1 bnd st2}
|g_1(\hat\al)|\le c\rho(m)\int_{\substack{|\tth|\le a\e^{-1/2}\\ |\CA(\tth)-\hal|\ge \de' \hat\al^{1/2}} }(1+(\CA(\tth)-\hal)^2)^{-1}\dd\tilde\theta.
\ee
Again by Lemma~\ref{lem:aux props*}, on the sets $\tth>0$ and $\tth<0$ we can change variables $\tth\to r=\CA(\tth)$, where $d\tth/dr \asymp \pm|r|^{-1/2}$, to obtain
\be\label{g1 bnd st3}
|g_1(\hat\al)|\le c\rho(m)\int_{\substack{r>0\\|r-\hal|\ge \de' \hat\al^{1/2}} }(1+(r-\hal)^2)^{-1}r^{-1/2}\dd r.
\ee
By an easy calculation,
%
\be\label{model ints}
\int_{\hal+\de' \hal^{1/2}}^\infty\frac{\dd r}{(r-\hal)^2 r^{1/2}},\ 
\int_0^{\hal-\de'\hat\al^{1/2}}\frac{\dd r}{(\hal-r)^2 r^{1/2}}=O(\hat\al^{-1}),
\ee
and the lemma is proven.

\subsection{Proof of Lemma \ref{lem: nozero}}
By \eqref{*r props II}, we need to establish that $|\e^{1/2}\hat\al/\tth|\le c$ for some sufficiently small $c>0$. 
Given that $\tth$, $\hal$, satisfy the conditions in \eqref{tth equiv}, \eqref{*r props I} and \eqref{*r props II} imply
\be
\left|\frac{\e^{1/2}\hat\al}{\tth}\right|
\le \frac{\e^{1/2}\hat\al}{\min_{|R(\tilde\theta,\hat\al)|< \de\hat\al^{1/2}}|\tth|}
\le c_1\frac{\e^{1/2}\hat\al}{(\hat\al-c_2\de\hat\al^{1/2})^{1/2}}
\ee
for some $c_{1,2}>0$. Squaring both sides gives
\be
\frac{\e\hat\al^2}{\tth^2}\le c_1^2\frac{\e\hat\al^2}{\hat\al-c_2\de\hat\al^{1/2}}=c_1^2\frac{\al}{1-c_2(\de/\hat\al^{1/2})}.
\ee
Given that $\al\in\Omega$, $\Omega$ can be made as small as we like (by selecting $a>0$ small), $\hat\al\ge\hal_*$, and $\hal_*$ can be made as large as we like, \eqref{tth equiv} is proven.

\subsection{Proof of Lemma \ref{g4 bnds}} 
We will consider only $g_4^+$ (i.e., $n>0$), since estimating $g_4^-$ is completely analogous. For simplicity, the superscript $+$ is omitted.

By Lemmas~\ref{lem:aux props*} and \ref{lem: nozero}, $\pa_{\tth}R\asymp\tth$ and $\CA^{-1}(\hal)\asymp\hal^{1/2}$, so the coefficient in front of the integral in \eqref{g4pl r} is $O(\hal^{-1/2})$. Suppose first that one of the intervals $R_{n_0}=[r_{2n_0},r_{2n_0+1}]$ contains zero. Then $r_{2n_0}\le 0\le r_{2n_0+1}$ and  
\be\label{psi alt int}
J:=\int_{r_{2n_0}}^{r_{2n_0+1}} \tilde\psi_m(r+h)\dd r=-\left(\int_{-\infty}^{r_{2n_0}}+\int_{r_{2n_0+1}}^{\infty}\right) \tilde\psi_m(r+h)\dd r.
\ee
Therefore, by \eqref{four-coef-bnd},
\be\label{J rh bound}
J\le \frac{c\rho(m)}{1+\min(|r_{2n_0}+h|,|r_{2n_0+1}+h|)}.
\ee
If either $r_{2n_0}<-\de\hal^{1/2}$ or $r_{2n_0+1}>\de\hal^{1/2}$, then the corresponding limit is replaced by either $-\de\hal^{1/2}$ or $\de\hal^{1/2}$, as needed. If both $r_{2n_0}$ and $r_{2n_0+1}$ exceed the limits, then $J=O(\hal^{-1/2})$. By Lemma~\ref{lem:aux props*}, $r_n\asymp v_n-\hal$, $v_n=\CA(u_n)$, and \eqref{g2 asymp} is proven.

Recall that $u_n>0$ if $n>0$. Then $u_n\asymp n$ and, by Lemma~\ref{lem:aux props*}, $v_n\asymp n^2$. This implies that, on average, the distance between consecutive $v_n$ increases as $n\to\infty$. In turn, this means that $v_n-\hal$ stays bounded for a progressively smaller fraction of $\hal$ as $\hal\to\infty$. Since the term $h$ is uniformly bounded, it can be omitted from \eqref{J rh bound} to better reflect the essence of the estimate. 

Contribution of all remaining intervals located on one side of zero $[r_{2n},r_{2n+1}]\subset (r_{2n_0+1},\de\hal^{1/2}]$, $n>n_0$, and $[r_{2n},r_{2n+1}]\subset [-\de\hal^{1/2},r_{2n_0})$, $0\le n<n_0$, can be estimated in a similar fashion:
\be\label{psi alt pos}
\sum_{n>n_0}\int_{r_{2n}}^{r_{2n+1}} |\tilde\psi_m(r+h)|\dd r
\le \int_{r_{2n_0+1}}^{\infty} |\tilde\psi_m(r+h)|\dd r\le \frac{c\rho(m)}{1+r_{2n_0+1}},
\ee
and the same way for the other set of intervals. This proves \eqref{g2 asymp}. 

%

To prove \eqref{g2derasymp}, we first collect some useful results, which follow from Lemma~\ref{lem:aux props*}:
\be\label{useful res}
\pa_{\hat\al}\CA^{-1}(\hat\al)=O(\hat\al^{-1/2}),\ 
\pa_{\hat\al} r_n=O(1).
\ee
Differentiating $g_4(\hat\al)$ in \eqref{g4pl r} and using \eqref{*r props II}, \eqref{*hF est}, and \eqref{useful res} gives
\be\label{g4 deriv}\begin{split}
&\pa_{\hat \al}g_4^\pm(\hat\al)\\
&\le 
c\rho(m)\left(\left(\frac{\e^{1/2}}{\hat\al}+\frac{1}{\hat\al^{3/2}}\right)\frac1{1+r_{\text{min}}(\hal)}
+\frac{1}{\hat\al^{1/2}}\frac1{1+r_{\text{min}}^2(\hal)}+\frac{\e^{1/2}}{\hal}+\frac{\e}{\hal^{1/2}}\right),
\end{split}
\ee
and the desired result follows by keeping only the dominant terms. Here we have used that $\e\hat\al=O(1)$, and there are finitely many $n$ such that $[r_{2n},r_{2n+1}]$ intersects the set $|r|\le\de\hat\al^{1/2}$. The last claim is proven by finding all $u_n\asymp n^2$ that satisfy $\hal-\de'\hal^{1/2}\le\CA(u_n)\le\hal+\de'\hal^{1/2}$ for some $\de'>0$.

\subsection{Proof of Lemma \ref{lem: on Vn}} 
We consider only $g_4^+$, the proof for $g_4^-$ is analogous. 
By Lemma~\ref{lem: Vn props},
$\hat\al\asymp v_n$ if $\hat\al\in V_n$. Also, at the right endpoint of $V_n$, by \eqref{g2 asymp}
\be\label{g4 endpt prf}
|g_4^+(v_{n+1/2})| \le c \rho(m) \left(v_n^{-1/2}(1+(v_{n+1}-v_n))^{-1}+v_n^{-1}\right)=\rho(m)O(n^{-2}).
\ee
By Lemma~\ref{g4 bnds}, with $\hat\al=v_n+\tau$, $\hat\al\in V_n$,
\begin{align}
\label{g4 bnd2}
&|g_4^+(\hat\al)| \le c \rho(m) \left(v_n^{-1/2}(1+\tau)^{-1}+v_n^{-1}\right),\\
\label{g2der bnd2}
&|\pa_{\hat\al}g_4^+(\hat\al)| \le c \rho(m) \left(v_n^{-1/2}(1+\tau^2)^{-1}
+\e^{1/2}v_n^{-1}\right)
\end{align}
if $v_{n-1/2}\ge \hat\al_*$. Then
\be\label{g4 ints prf}\begin{split}
&\int_{V_n}|g_4^+(\hat\al)|\dd\hat\al \le c \rho(m) \left(v_n^{-1/2}\int_0^{O(n)}\frac{\dd\tau}{1+\tau}+v_n^{-1}O(n)\right)
=\rho(m)O(n^{-1}\ln n),\\
&\int_{V_n}|\pa_{\hat\al}g_4^+(\hat\al)|\dd\hat\al \le c \rho(m) \left(v_n^{-1/2}O(1)
+\e^{1/2}v_n^{-1} O(n)\right)
=\rho(m)O(n^{-1}),
\end{split}
\ee
because $n=O(\e^{-1/2})$, and the lemma is proven.

\section{Analysis of exceptional cases}\label{sec: spec cases}
Consider now possible violations of the inequalities $0<\al_l<\al_*<\al_r\le \pi/2$ (see at the beginning of sections~\ref{ssec:bulk} and \ref{ssec: al*}). A violation may happen only for finitely many $m$. If $\al_l<0$ (i.e., the interval $(\phi(0),\phi(\al_*))$ does not contain an integer), we consider the interval $[0,\al_*]$ instead of $[\al_l,\al_*]$. The analogues of \eqref{I*l sub v0}--\eqref{thpr*} become
\be\label{I*l sub}\begin{split}
&N=O(\tal_*)=O(\e^{-1/2});\
[0,\tal_*]=\cup_{n=0}^{N-1} [nL,(n+1)L],\ \tal_*=NL;\\
&\langle \vartheta'(\hal)\rangle\ge \min\left(\langle m\kappa r_x\rangle,\langle m\mu\rangle\right),\
\al\in[0,\al_*].
\end{split}
\ee
The analogue of \eqref{Wjm est stp1} is (cf. \eqref{G bnd comb})
\be\label{Wjm est special}\begin{split}
W_m([0,\al_*])\le &c\frac{\e^{1/2}\rho(m)}{\min\left(\langle m\kappa r_x\rangle,\langle m\mu\rangle\right)}\sum_{n=0}^{N-1}\frac1{1+(N-(n+1))L}=O(\e^{1/2}\ln(1/\e)).
\end{split}
\ee
Here we have used that the number of different values of $m$ for which $\al_l<0$ and the above estimate applies is finite. 

If $\al_r>\pi/2$ (which implies $\langle \phi(\pi/2)\rangle>0$), we consider the interval $[\al_*,\pi/2]$ instead of $[\al_*,\al_r]$. The analogue of \eqref{I*l sub} becomes
\be\label{I*r sub}\begin{split}
&N=O(\e^{-1/2}),\\ 
&[\tal_*,\e^{-1/2}\pi/2]=\cup_{n=0}^{N-1} [\tal_*+nL,\tal_*+(n+1)L],\ \tal_*+NL=\e^{-1/2}\pi/2,\\
&\langle \vartheta'(\hal)\rangle\ge \min\left(\langle m\kappa r_x\rangle,\langle \phi(\pi/2)\rangle\right),\
\al\in[\al_*,\pi/2].
\end{split}
\ee
The analogue of \eqref{Wjm est special} is
\be\label{Wjm est sp-r}
W_m([\al_*,\pi/2])\le c\frac{\e^{1/2}\rho(m)}{\min\left(\langle m\kappa r_x\rangle,\langle \phi(\pi/2)\rangle\right)}=O(\e^{1/2}).
\ee
Here we used that $\tal=O(\e^{-1/2})$ if $\al\ge\al_*$.

This completes the analysis of the case $\al_*<\pi/2$. If $\al_*>\pi/2$, then $\al_l>\pi/2$ when $m\gg1$ is sufficiently large. Hence the only relevant case is when $\al_l<\pi/2$. As before, this happens for finitely many $m$. If $\al_l>0$, the interval we should consider is $[\al_l,\pi/2]$. In this case, the estimates in \eqref{Wjm est stp1}, \eqref{Wjm est v2} still hold (some of the terms in the sum are not necessary). Given that $m$ is bounded, the estimates imply $W_m([\al_l,\pi/2])=O(\e^{1/2}\ln(1/\e))$. If $\al_l<0$, the relevant interval is $[0,\pi/2]$. Arguing similarly to \eqref{I*l sub}, \eqref{Wjm est special}, it is obvious that $W_m([0,\pi/2])=O(\e^{1/2}\ln(1/\e))$.

\section{Proof of Lemma~\ref{lem:partial res}}\label{sec new ker prf}

Pick any $x$ sufficiently close to $x_0$. All the estimates below are uniform with respect to $x$ in a small (but fixed) size neighborhood, so the $x$-dependence of various quantities is omitted from notation. 

Let $\Omega$ be the set of all $\al\in[-\pi/2,\pi/2]$ such that the lines $\{y\in\br^2:\,(y-x)\cdot\vec\al=0\}$ intersect $\s$. Let $\theta=\Theta(\al)$, $\al\in\Omega$, be determined by solving $(y(\theta)-x)\cdot\vec\al=0$. By using a partition of unity, if necessary, we can assume that the solution is unique.  By assumption, the intersection is transverse for any $\al\in\Omega$ (up to the endpoints). Hence 
$|\Theta'(\al)|=|y(\Theta(\al))-x|/|\vec\al\cdot y'(\Theta(\al))|$ and 
\be\label{theta deriv min}
|\vec\al\cdot y'(\Theta(\al))|,|\Theta'(\al)|\asymp 1,\ \al\in\Omega.
\ee

Transform the expression for $A_m$ (cf. \eqref{recon-ker-v2}) similarly to \eqref{A-simpl}:
\be\label{A-simpl-app}
\begin{split}
A_m(\al,\e)=\frac1\e\int_{-a}^a\int_0^{\e^{-1}H_\e(\theta)}\tilde\psi_m\left(\frac{\vec\al\cdot (y(\theta)-x)}\e+\hat t\cos(\theta-\al)\right)F(\theta,\e \hat t)\dd\hat t\dd\theta,
\end{split}
\ee
where $F$ is the same as in \eqref{A-simpl}. Setting $\tilde\theta=(\theta-\Theta(\al))/\e^{1/2}$,
\eqref{A-simpl-app} becomes
\be\label{A-simpl_2}
\begin{split}
A_m(\al,\e)=&\e^{-1/2}\int\int_0^{\e^{-1}H_\e(\theta)}\tilde\psi_m\left(\frac{\vec\al\cdot (y(\theta)-y(\Theta(\al)))}\e+\hat t\cos(\theta-\al)\right)\\
&\times F(\theta,\e \hat t)\dd\hat t\dd\tilde \theta,\quad
\theta=\Theta(\al)+\e^{1/2}\tilde\theta,\ \al\in\Omega.
\end{split}
\ee
Note that in this appendix the relation between $\tth$ and $\theta$ is different from that in the rest of the paper (see \eqref{scales}).

Due to \eqref{four-coef-bnd}, we can integrate with respect to $\tilde\theta$ over any fixed neighborhood of $0$:
\be\label{A-simpl-st2}
\begin{split}
&A_m(\al,\e)\\
&=\e^{-1/2}\int_{-\de}^{\de}\int_0^{\e^{-1}H_\e(\theta)}\tilde\psi_m\left(\frac{\vec\al\cdot y'(\Theta(\al))}{\e^{1/2}}\tilde\theta+O(\tilde\theta^2)+\hat t\cos(\Theta(\al)-\al)+O(\e^{1/2})\right)\\
&\quad\times \left(F(\Theta(\al),0)+O(\e^{1/2})\right)\dd\hat t\dd\tilde \theta+\rho(m)O(\e^{1/2})
,\quad \theta=\Theta(\al)+\e^{1/2}\tilde\theta,\ \al\in\Omega,
\end{split}
\ee
for some $\de>0$ sufficiently small. Using \eqref{four-coef-bnd} it is easy to see that the terms $O(\e^{1/2})$ and $O(\tilde\theta^2)$ can be omitted from the argument of $\tilde\psi_m$ without changing the error term:
\be\label{A-simpl-st3}
\begin{split}
A_m(\al,\e)
=&\frac{F(\Theta(\al),0)}{\e^{1/2}}\int_{-\de}^{\de}\int_0^{H_0(\e^{-1/2}\Theta(\al)+\tilde\theta)}\tilde\psi_m\left(a(\al)\e^{-1/2}\tilde\theta+b(\al)\hat t\right)\dd\hat t\dd\tilde \theta\\
&+\rho(m)O(\e^{1/2}),\quad 
a(\al):=\vec\al\cdot y'(\Theta(\al)),\ b(\al):=\cos(\Theta(\al)-\al),\ \al\in\Omega.
\end{split}
\ee
By \eqref{theta deriv min}, $a(\al)$ is bounded away from zero on $\Omega$.  By the last equation in \eqref{psi-props-0}, $\int\tilde\psi_m(\hat t)\dd\hat t=0$ for all $m$, so we can replace the lower limit in \eqref{A-simpl-st3} with any value independent of $\tilde\theta$. Again, we use here that the contribution to the integral with respect to $\tilde\theta$ over $\br\setminus(-\de,\de)$ is $\rho(m)O(\e^{1/2})$. We choose the lower limit to be $H_0(\e^{-1/2}\Theta(\al))$:
\be\label{A-simpl-st4}
\begin{split}
A_m(\al,\e)
=&\frac{F(\Theta(\al),0)}{\e^{1/2}}\int_{-\de}^{\de}\int_{H_0(\e^{-1/2}\Theta(\al))}^{H_0(\e^{-1/2}\Theta(\al)+\tilde\theta)}\tilde\psi_m\left(a(\al)\e^{-1/2}\tilde\theta+b(\al)\hat t\right)\dd\hat t\dd\tilde \theta\\
&+\rho(m)O(\e^{1/2}),\ \al\in\Omega.
\end{split}
\ee
Neglecting the big-$O$ term in \eqref{A-simpl-st4} leads to a term of magnitude $O(\e^{1/2})$ in $\frec$. Clearly,
\be\label{A-outside}
A_m(\al,\e)=\rho(m)O(\e),\ \al\in[-\pi/2,\pi/2]\setminus\Omega.
\ee
%

Define
\be\label{sum-k}
\begin{split}
B_m(\tilde\theta,\e)
:=&\Delta\al\sum_{\al_k\in\Omega} e\left(-m q_k\right) F(\Theta(\al_k),0)\\
&\qquad\times\int_{H_0(s_k)}^{H_0(s_k+\tilde\theta)}\tilde\psi_m\left(a(\al_k)\e^{-1/2}\tilde\theta+b(\al_k)\hat t\right)\dd\hat t,\\ 
s_k:=&\e^{-1/2}\Theta(\al_k),\ q_k=\vec\al_k\cdot x/\e.
\end{split}
\ee
From \eqref{recon-ker-v2} and \eqref{A-simpl-st4}--\eqref{sum-k},
\be\label{bm to frec}
\frec(x)=c\e^{-1/2}\sum_m e(m\bar p/\e) \int_{-\de}^{\de} B_m(\tilde\theta,\e)\dd\tth+O(\e^{1/2}).
\ee
Clearly,
\be\label{Bm est}
\begin{split}
|B_m(\tilde\theta,\e)|
\le &c \Delta\al\sum_{\al_k\in\Omega} \left|\int_{H_0(s_k)}^{H_0(s_k+\tilde\theta)}\tilde\psi_m\left(a(\al_k)\e^{-1/2}\tilde\theta+b(\al_k)\hat t\right)\dd\hat t\right|.
\end{split}
\ee
Define similarly to \eqref{chi fn}:
\be\label{H-def alt}
\chi_{t_1,t_2}(t):=\begin{cases}1,&t_1\le t\le t_2 \text{ or }t_2\le t\le t_1,\\
0,&\text{ otherwise},
\end{cases}
\ee
and rewrite \eqref{Bm est} in the form
\be\label{sum-k alt}
\begin{split}
|B_m(\tilde\theta,\e)|
\le &c\int \Delta\al\sum_{\al_k\in\Omega}\left|\tilde\psi_m\left(a(\al_k)\e^{-1/2}\tilde\theta+b(\al_k)\hat t\right)\right| \chi_{H_0(s_k),H_0(s_k+\tilde\theta)}(\hat t)\dd\hat t\\
=&\int g_{\e,m}(\tilde\theta,\hat t) \dd\hat t,
\end{split}
\ee
where
\be\label{G and a_k b_k}
\begin{split}
&g_{\e,m}(\tilde\theta,\hat t):=
\Delta\al \sideset{}{'}\sum_k \left|\tilde\psi_m\left(a(\al_k)\e^{-1/2}\tilde\theta+b(\al_k)\hat t\right)\right|.
\end{split}
\ee
The prime next to the summation sign means that the sum is over all $\al_k\in\Omega$ such that either $H_0(s_k)< \hat t<H_0(s_k+\tilde\theta)$ or $H_0(s_k+\tilde\theta) < \hat t<H_0(s_k)$
By ignoring a set of measure zero, we assume here and below that $\hat t\not= H_0(s_k)$ and $\hat t\not=H_0(s_k+\tilde\theta)$ for any $k$. 

Our argument implies that index $k$ appears in the sum in \eqref{G and a_k b_k} only if $s_k$ and $s_k+\tth$ are on different sides of at least one $u_n$, e.g. $s_k\in U_n$ and $s_k+\tth\not\in U_n$ or vice versa (see the text following \eqref{U intervals}). By assumption $H3$, the number of $u_n$ on an interval of length $O(\e^{-1/2})$ is $O(\e^{-1/2})$ uniformly in $\hat t$. Hence
\be\label{dbl sum}
\begin{split}
&g_{\e,m}(\tilde\theta,\hat t)\le
\Delta\al
\sum_{|u_n|\le O(\e^{-1/2})} \sum_{|s_k-u_n|\le |\tth|}
\left|\tilde\psi_m\left(a(\al_k)\e^{-1/2}\tilde\theta+b(\al_k)\hat t\right)\right|.
\end{split}
\ee
Fix any $n$. By \eqref{theta deriv min}, there are no more than $1+O(\e^{-1/2}|\tilde\theta|)$ values of $k$ such that $|s_k-u_n|\le \tth$.
Now we estimate $g_{\e,m}$ using \eqref{four-coef-bnd}, \eqref{theta deriv min}, \eqref{G and a_k b_k}, and \eqref{dbl sum}:
\be\label{W est}
\begin{split}
&g_{\e,m}(\tilde\theta,\hat t)\le
c\rho(m)\e^{1/2}\frac{\max(1,\e^{-1/2}|\tth|)}{1+(\tilde\theta^2/\e)}.
\end{split}
\ee
Since the set of $\hat t$ values (the range of $H_0$) is uniformly bounded, the same estimate as in \eqref{W est} applies to $B_m(\tilde\theta,\e)$ (cf. \eqref{sum-k alt}). Using \eqref{bm to frec}, we finish the proof:
\be\label{last bnd}
\sum_m \e^{-1/2}\int_{-\de}^{\de} |B_m(\tilde\theta,\e)|\dd \tilde\theta
=O(\e^{1/2}\ln(1/\e)).
\ee

%

\bibliographystyle{plain}
\bibliography{My_Collection}
\end{document}